\RequirePackage{fix-cm}

\documentclass[smallextended,final]{svjour3}       

\smartqed  
\usepackage{mathptmx}
\usepackage{graphicx}
\usepackage{natbib}
\usepackage{longtable}
\usepackage{amsmath}
\usepackage{tabularx}
\usepackage{booktabs}
\usepackage{hhline}
\usepackage{multirow}
\usepackage{algorithm, algpseudocode, color, soul}
\usepackage{versions}
\usepackage[section]{placeins}
\usepackage{mathpazo}
\usepackage{enumerate}
\usepackage{upgreek}
\usepackage{amssymb}
\usepackage{subfigure}
\usepackage{hyperref}

\newcommand{\tagarray}{%
\mbox{}\refstepcounter{equation}%
$(\theequation)$%
}
\newcolumntype{L}[1]{>{\raggedright\arraybackslash}p{#1}} 
\newcolumntype{C}[1]{>{\centering\arraybackslash}p{#1}} 
\newcolumntype{R}[1]{>{\raggedleft\arraybackslash}p{#1}} 

\begin{document}

\title{Stochastic project management: Multiple projects with multi-skilled human resources
}


\author{Thomas Felberbauer         \and
        Walter J. Gutjahr \and Karl F. Doerner
}


\institute{T. Felberbauer \at
              Department of Media and Digital Technologies, St. P\"olten University
			  of Applied Sciences, St. P\"olten, Austria\\
              \email{thomas.felberbauer@fhstp.ac.at}           
			\and  W. J. Gutjahr \at
              Department of Statistics and Operations Research, University of Vienna, Vienna, Austria\\
              \email{walter.gutjahr@univie.ac.at}
           \and
            K. F. Doerner \at
              Department of Business Decisions and Analytics, Data Science @ University of Vienna, Vienna, Austria\\
              \email{karl.doerner@univie.ac.at}
}

\date{Received: date / Accepted: date}

\maketitle

\begin{abstract}
This paper presents two stochastic optimization approaches for simultaneous project scheduling and personnel planning, extending a deterministic model previously developed by Heimerl and Kolisch. For the problem of assigning work packages to multi-skilled human resources with heterogeneous skills, the uncertainty on work package processing times is addressed. In the case where the required capacity exceeds the available capacity of internal resources, external human resources are used. The objective is to minimize the expected external costs. The first solution approach is a ``matheuristic'' based on a decomposition of the problem into a project scheduling subproblem and a staffing subproblem. An iterated local search procedure determines the project schedules, while the staffing subproblem is solved by means of the Frank-Wolfe algorithm for convex optimization. The second solution approach is Sample Average Approximation where, based on sampled scenarios, the deterministic equivalent problem is solved through mixed integer
programming. Experimental results for synthetically generated test instances inspired by a real-world situation are provided, and some managerial insights are derived.

\end{abstract}

\section{Introduction}
\label{intro}
In project management, the two important tasks of project scheduling on the one hand and of personnel planning on the other hand are usually not faced separately from each other, but rather in a simultaneous or interleaved manner.
This interaction considerably increases the complexity of the planning process, with the consequence that project managers frequently wish to get computational support alleviating the cognitive burden of the combined scheduling/staffing decision. Decision support tools addressing the combined problem have been proposed in the literature, but typically they use deterministic models. Despite the fact that managers and scientists recognize the importance of uncertainty within project management, few works take this aspect into account when proposing methods for quantitative decision support for the combined problem described above, obviously for the reason that, as argued, this problem is already hard in the deterministic case.
In this study we highlight the potential of stochastic optimization applied to combined project scheduling and staffing, and compared two different methods of solving the proposed model computationally.

\begin{sloppypar}
The project scheduling and staffing model proposed by~\cite{Heimerl2010a} deals with the problem of assigning multi-skilled human resources to work, while taking into account resource-specific and heterogeneous skill efficiencies. The objective is to minimize the costs for internal and external personnel. \cite{Heimerl2010a} provide a mixed-integer programming (MIP) formulation with a tight LP bound and compare its performance, using the MIP solver of CPLEX, to that of simple heuristics. Additionally, they investigate the influence of different parameters, such as the time window size, the utilization, and the number of skills per resource, on the personnel costs. In \cite{Kolisch2012} the authors extend their previously presented model and develop a hybrid metaheuristic as an innovative solution method. The proposed solution method separates the problem into a scheduling and a staffing subproblem, where the staffing problem is solved through a generalized network simplex algorithm. Results show that the hybrid metaheuristic outperforms the MIP solver. In \citet{Felberbauer2016b}, the models from~\cite{Heimerl2010a} and~\cite{Kolisch2012} are extended by considering labor contracts. For the problem of project scheduling and staffing with labor contracts, \citet{Felberbauer2016b} propose a hybrid metaheuristic combining iterated local search and a greedy staffing heuristic.
\end{sloppypar}

\cite{Bergh2013} observe that most papers on personnel scheduling problems still appear to feature a deterministic approach and advise researchers to consider stochastic problem versions. They highlight the importance of uncertainty, explicitly mentioning volatile demand, last-minute changes, and rescheduling based on new information, as interesting new research topics.
A recent article on workforce planning by \cite{Bruecker2015} states that the small number of papers in this field taking uncertainty into account is alarming. They refer to the general consensus that uncertainty is ubiquitous in real workforce planning, but uncertainty still remains relegated to the ``future research'' section of many papers. The authors mention the high complexity of integrating uncertainty within optimization models and techniques as one reason for its inadequate representation in the literature.

Considering uncertainty in project management and its subordinate planning levels, a recent work of \citet{Barz2014} describes resource assignments in the telecommunication industry. The authors investigate a hierarchical, multi-skill resource assignment problem and use a discrete Markov decision process model for incoming jobs over an infinite time horizon.
\citet{Gutjahr2013} present a stochastic optimization model for project selection and project staffing. They assume that both the returns and the required efforts of the selected projects are random variables. The problem is decomposed into a project selection problem and a staffing subproblem. For the computational solution of the two problems, an adapted version of Variable Neighborhood Search and a Frank-Wolfe type algorithm, respectively, are used.
\citet{Artigues2013} propose a robust optimization approach to resource-constrained project scheduling with uncertain activity duration, assuming that the decision-maker cannot associate probabilities with possible activity durations. The authors describe how robust optimization can be applied to project scheduling under uncertainty, and they develop a scenario-relaxation algorithm and a heuristic solution method that solves medium-sized instances.
\citet{Gutjahr2015} proposes a model for stochastic multi-mode resource-constrained project scheduling under risk
aversion with the two objectives of makespan and costs. Activity durations and costs are modelled as random variables.
For the scheduling part of the decision problem, the class of early-start
policies is considered. A further decision to be made concerns the assignment of execution modes to activities.
To take risk aversion into account, an approach of optimization under multivariate stochastic
dominance constraints is adopted. For the resulting bi-objective stochastic
integer programming problem, the Pareto frontier is determined by means of an exact solution method,
incorporating a branch-and-bound technique.

\citet{Ingels2017employee,ingels2017} investigate the assignment of employees to cover the staffing requirements for specific skills and shifts. Their research focuses on  improving the short-term adjustment capability of the shift roster by maximising
the substitutability of employees as a means to increase the robustness of a project plan. The authors propose a three-step methodology including a two-phase pre-emptive programming approach. Both uncertainty of demand and uncertainty of capacity are considered by stochastic models. A personnel shift roster for a medium-term period offering sufficient flexibility for between-skill substitution and within-skill substitution is provided.

\begin{sloppypar}
The articles~\citet{Heimerl2010a}, \citet{Kolisch2012}, \citet{Felberbauer2016b} and \citet{Gutjahr2013} constitute the starting point for the present paper. In line with these articles, we assume that in the case where the work time demand exceeds the available capacity, external capacity is used, e.g. by hiring external personnel. However, we introduce the following new features:
(a) Contrary to the deterministic models presented in~\citet{Heimerl2010a, Kolisch2012, Felberbauer2016b}, we model the required work time demand as {\em stochastic}. On the other hand, in order to keep our model compact, we simplify the model of~\citet{Heimerl2010a, Kolisch2012} by neglecting the possibility of overtime work.
(b) Contrary to \citet{Gutjahr2013}, we address the {\em project scheduling} decision. This issue introduces a relevant additional source of computational complexity into the problem, since project scheduling problems (already in a simple, deterministic context) are notoriously hard. In particular, the demand information does not appear here per project, but in a more fine-grained way for each activity of a project.  On the other hand, compared to \citet{Gutjahr2013}, we do not deal here with the project selection aspect.
\end{sloppypar}

\begin{sloppypar}
To cope with the computational challenge of the proposed model, we develop two solution approaches and compare their performance. The first, ``matheuristic'' approach uses a metaheuristic for the scheduling part of the problem as well as an exact solution procedure inspired by~\citet{Gutjahr2013} for the staffing part. It is noteworthy that the staffing part constitutes a nonlinear (though convex) optimization problem. The second approach follows a completely different strategy by using Sample Average Approximation to achieve a linearization of the model, and by solving the resulting ``deterministic equivalent'' problem, which remains a combined scheduling-staffing problem, with the aid of a MIP solver. A comparison of these two ways to address a mixed-integer stochastic optimization problem with a convex ``lower level'' problem may also be interesting for applications outside the area of project management.
\end{sloppypar}

Substantial literature focuses on the deterministic case of project scheduling and staffing. However, to the best of our knowledge, no previous studies have addressed the scheduling and staffing of multiple projects with multi-skilled resources holding heterogeneous skill efficiencies, under consideration of uncertain work time demand. The contribution of our article is:
\begin{enumerate}
\item We adapt the model of \cite{Kolisch2012} by considering stochastic demand.
\item We present two solution procedures, namely Sample Average Approximation and a matheuristic, and compare their performance.
\item We analyze deterministic and stochastic planning and discuss cost estimation accuracy and the value of the stochastic solution.
\item We formulate managerial insights from the stochastic planning approach.
\end{enumerate}

The paper is organized as follows: In Section~\ref{sec:problem} the stochastic optimization model is introduced. Section~\ref{sec:solutionprocedure} presents the developed solution procedures. The structure of our test instances and the experimental design are described and the results of the computational experiments are discussed in Section~\ref{sec:experiments}. Finally, in Section~\ref{sec:Conclusion and further research} we formulate some managerial insights that can be drawn from the study of the stochastic optimization problem, and we mention some interesting topics for future research.

\section{Problem formulation}
\label{sec:problem}
In this section the stochastic project scheduling and staffing problem is presented.

\textbf{\textit{Projects:}}
$P$ is a set of {\em projects} that must be completed during the planning horizon $T$ (e.g., one year). All projects are independent of each other but compete for the same resources. The planning horizon $T$ comprises discrete time periods $t$ (e.g., months). Each project $p$ $(p \in P)$ has a work schedule defined by a sequence of consecutive {\em activities} $q=1,\ldots,d_p$, where an activity is defined as that part of a project that is executed during one time period $t$. In this way, we discretize the overall project into consecutive pieces, each of the same length; as it will be discussed below, interruptions will be allowed. During activity~$q$, project $p$ demands $D_{psq}$ units of work requiring {\em skill} $s$ $(s \in S)$, where $S$ is a set of skills. The amount of work in a skill~$s$ performed in a certain activity~$q$ of a certain project~$p$ is called a {\em work package}. The sizes~$D_{psq}$ of the work packages, which we can also interpret as work contents or efforts, are considered as random variables. This reflects the very frequently occurring situation that the required efforts~$D_{psq}$ are unknown at the beginning of the planning process. Their true realizations become only known during project execution; we assume, however, that the distribution of the variables $D_{psq}$ is known in advance(or can at least be estimated): it is described by the probability density function (PDF) $h_{psq}$ $(p \in P, s \in S, 1 \le q \le d_p)$. By using a PDF, we suppose that the distribution of~$D_{psq}$ is a continuous distribution, such as a triangular or a beta distribution.

Moreover, we assume that a {\em time window} $[ES_p, LS_p]$ for the start of each project $p \in P$ is pre-defined. Therein, $ES_p$ and $LS_p$ denote the index~$t$ of the earliest time period and of the latest time period, respectively, in which project~$p$ can be started. The {\em latest finish period} of project $p$, denoted $LF_p$, is defined as $LF_p=LS_p+d_p-1$. The {\em time window size} $\gamma$ of project $p$ is defined as $\gamma = LS_p - ES_p$. The earliest start time of activity~$q$ is denoted $ES_{pq}$ and defined as $ES_{pq}=ES_p+q-1$. The latest start time of activity~$q$ of project~$p$, denoted $LS_{pq}$, is equal to the latest finish time $LF_{pq}$ of this activity and is defined as $LS_{pq}=LF_{pq}=LS_p+q-1$. The activities $q=1, \dotsc ,d_p$ have to be processed in an ascending order. However, it is possible to {\em interrupt} the project between the discrete activities $q$ once the project has started. Nevertheless, the order of the activities as well as the start time and finish time constraints have to be respected. For example, if $ES_p=1$, $LS_p=3$ and $d_p=3$ for a given project~$p$, then it would be allowed to schedule the three activities of project~$p$ in time periods~1, 2, 3 (earliest possible schedule); some other feasible alternatives would be to schedule them in time periods 1, 4, 5, in time periods 2, 3, 5, or in time periods 3, 4, 5 (latest possible schedule). The schedule 2, 4, 6 would not be allowed, because it exceeds the latest finish period $LF_p = 5$.

\textbf{\textit{Resources:}}
Our model distinguishes between {\em internal} and {\em external} human resources. Each internal resource~$k$ taken from the set of all internal human resources~$K$ holds a subset of skills $S_k \subseteq S$. Conversely, from the perspective of a skill~$s$, the subset $K_s = \{k \in K \: | \: s \in S_k\}$ of~$K$ contains all resources that can perform skill $s$. If $s \in S_k$, we take account of different degrees in the performance of resource~$k$ in skill~$s$, specified by the {\em efficiency}~$\eta_{sk}$ of resource~$k$ in skill~$s$. The higher the efficiency value~$\eta_{sk}$, the faster is resource~$k$ in executing a work package $D_{psq}$ requiring skill~$s$.
It is assumed that $\eta_{sk}$ is known at the beginning of the time horizon (the time of the decision) and does not change until its end. (In particular, we disregard effects as learning or knowledge depreciation). The values $\eta_{sk}$ are strictly positive for all $s \in S_k$. Whenever convenient, we extend the definition of the values~$\eta_{sk}$ to all pairs $(s,k)$ and set then $\eta_{sk} = 0$ for $s \notin S_k$. The quantity $\eta_{sk}$ represents a factor for the speed in which a certain skill is exerted. To distinguish between work content and actual time needed to perform a work package, the realization $d_{psq}$ of the random variable $D_{psq}$ will also be called \textit{effective work time}. To get the time actually needed by the considered resource, we have to divide the effective work time by $\eta_{sk}$, so the \textit{real work time} is $d_{psq}/\eta_{sk}$. For example, let us assume that the work package $d_{psq}$ requires 10 units of effective work time to be processed. (This would be the time a ``standard'' employee would need for it.) Then, if employee $k$ has an efficiency of $\eta_{sk}=1.25$ in skill~$s$, she or he needs only $d_{psq}/\eta_{sk}=10/1.25=8$ time units of \textit{real work time} to complete this work package.

\begin{sloppypar}
\textbf{\textit{Capacities:}}
Each (internal) resource~$k$ has a limited capacity of $a_{kt}$ during time period~$t$; these values are given in advance.
Similar to \citep{Gutjahr2013}, we assume that internal resources earn fixed wages per period. On our assumption that there is no overtime work, the costs for these wages sum up to a constant value, so they have no influence on the optimization problem. (Note that for the assumed employment type, even if an employee has idle times, the company has to pay for the hours where she or he is available and not only for the hours of scheduled work.) On the other hand, {\em external} resources are paid for the effective work they do. Our model assumes that they are available for each skill~$s$ to an unlimited extent at the expense of a cost rate~$c_s^e$ per unit of effective work time.

Since in our model, the work time demand $D_{psq}$ is stochastic, but the capacities of the internal resources are fixed and the project schedule has to be decided upon in advance as well, it is necessary to perform {\em recourse actions} after having observed the actual realizations of the $D_{psq}$ in order to be able to stick to the chosen schedule. We assume that whenever the internal capacity turns out as insufficient to cover the demand, the company resorts to external resources (experts in the required skills) and pays them at their given cost rates. The (expected) overall cost for the external resources has to be minimized.
\end{sloppypar}

\textbf{\textit{Decision variables:}}
For the stochastic optimization model, the decision variables $x_{ptsk} \geq 0$ define the amount of work performed by internal resource $k$ with skill $s$ in period $t$ for project $p$. Note that $x_{ptsk}$ is measured in effective work time and not in real work time. To define the times when the periods of a project are executed, the binary decision variables $z_{pqt}\in\{0,1\}$ are introduced: $z_{pqt} = 1$ if in time period~$t$, activity~$q$ of project~$p$ is executed, and 0 otherwise. Thus, the $x_{ptsk}$ variables define the staffing decision, whereas the $z_{pqt}$ variables define the scheduling decision.

\subsection{Stochastic optimization model}
\label{subsec:sm}

In (\ref{eq:objectivefunction1})~--~(\ref{eq:constraints6}) below, we present our stochastic optimization model in mathematical terms. For a project~$p$ and a time period~$t$, we let\\
\noindent\begin{tabularx}{\linewidth}{@{} L {80mm}L {20mm}>{\raggedleft\arraybackslash}X@{}}\\
\multicolumn{2}{c}{
$\displaystyle \uptau_{pt}= \lbrace q \in \lbrace 1,\dotsc,d_p\rbrace \: \vert \: ES_{pq} \leq t \leq LS_{pq}\rbrace$
} & \tagarray\label{eq:taupt}\\
&
\end{tabularx}
denote the set of all possible activities~$q$ of project $p$ that could lead to a resource demand in time period $t$.

\noindent\begin{tabularx}{\linewidth}{@{} L {80mm}L {20mm}>{\raggedleft\arraybackslash}X@{}}\\
\multicolumn{2}{c}{
\hspace*{-10ex}$\displaystyle\min_{z,x} \; \sum\limits_{p \in P} \: \sum\limits_{t = ES_p}^{LF_p} \: \sum\limits_{s \in S}  c_s^e \: E \bigg( \big[\sum\limits_{q \in \uptau_{pt}}D_{psq}z_{pqt}-\sum\limits_{k\in K_s} x_{ptsk}\big ]^+\bigg) \quad s.t.$
} & \tagarray\label{eq:objectivefunction1}\\
\\
$\displaystyle\sum\limits_{t=ES_{pq}}^{LS_{pq}} z_{pqt}=1$&
$\displaystyle \begin{aligned}
	p&\in P,\\
	q&=1,\dotsc,d_p
\end{aligned}$
&\tagarray\label{eq:constraints2}\\\addlinespace[2ex]
$\displaystyle\sum\limits_{t=ES_{pq}}^{LS_{pq}} t \, z_{pqt} <\displaystyle\sum\limits_{t=ES_{p,q+1}}^{LS_{p,q+1}} t \, z_{p,q+1,t}$&
$\displaystyle \begin{aligned}
	p&\in P,\\
	q&=1,\dotsc,d_p-1
\end{aligned}$
&\tagarray\label{eq:constraints3}\\\addlinespace[2ex]
$\displaystyle\sum\limits_{p \in P}\sum\limits_{s \in S} \frac{1}{\eta_{sk}}x_{ptsk} \leq a_{kt}$&
$\displaystyle\begin{aligned}
	k&\in K,\\
	t&=1, \dotsc, T
\end{aligned}$
&\tagarray\label{eq:constraints4}\\\addlinespace[2ex]
$\displaystyle x_{ptsk} \geq 0$ &
$\displaystyle \begin{aligned}
	p &\in P,\\
	t&=ES_p,\dotsc,LF_p,\\
	s &\in S,\\
	k &\in K_s
\end{aligned}$
&\tagarray\label{eq:constraints5}\\\addlinespace[2ex]
$\displaystyle z_{pqt}\in\{0,1\}$&
$\displaystyle \begin{aligned}
	p&\in P,\\
	q&=1,\dotsc,d_p,\\
	t&=ES_{pq},\dotsc,LS_{pq}
\end{aligned}$
&\tagarray\label{eq:constraints6}
\end{tabularx}
\vspace{1ex}

\begin{sloppypar}
The \textit{objective function} in Equation \eqref{eq:objectivefunction1} minimizes the expected external costs. The symbol $E$ denotes the mathematical expectation and $x^+$ stands for $\max(x,0)$. The external costs result from that part of the stochastic demand that is not covered by the internally scheduled capacity.
The constraints in Equation \eqref{eq:constraints2} ensure that activity $q$ of project $p$ is scheduled exactly once within the pre-defined time window $[ES_{pq},LS_{pq}]$.
The condition that activities $q$ must be processed in an ascending order is guaranteed by the constraints in Equation \eqref{eq:constraints3}.
Equation~\eqref{eq:constraints4} formulates the capacity constraints for each resource $k$ and time-period $t$: observe that for a resource with efficiency~$\eta_{sk}$, an effective work time of $x_{ptsk}$ entails a real work time of $x_{ptsk}/\eta_{sk}$.
Finally, the decision variables are defined in Equations~\eqref{eq:constraints5} and \eqref{eq:constraints6}.
\end{sloppypar}

Let us remark that in some cases, external costs grow faster than linearly in dependence of the outsourced work, for example because a selected supplier company has limited capacities itself. To generalize the model above to this situation, the given nonlinear external cost function could be approximated by a piecewise linear function, which can be dealt with by a simple extension of our Eq.~(\ref{eq:objectivefunction1}). Numerical solution techniques for this more complex model are a topic of future research.


As it is seen from the formulation above, our model can be viewed as a tactical two-stage stochastic optimization model with project plan and work assignment as first-stage decisions and outsourcing to external resources as the possible recourse actions. During the 'real' project execution, we assume that we get updated information about the realizations of the work package processing times. Our assumed recourse in the case of insufficient internal capacity is the hiring of external resources, where we assume that the information about the work package processing times arrives early enough to hire the necessary external resources. The presented tactical planning approach could be extended including more operative decisions by concerning other recourse actions such as re-assignment and re-scheduling of internal resources. This could lead to better solutions and new interesting managerial insights but would lead to a more sophisticate planning approach.

Finally, let us recall that our model assumes that the activities~$q$ of each project~$p$ have to be arranged in a pre-specified linear order related to \citet{Kolisch2012}. Nevertheless, due to the requirement of independent work package processing time distributions an additional assumption, where for each project~$p$ and in each time period~$t$ not more than one activity~$q$ of project~$p$ is allowed to be scheduled (no parallel activity processing per project), is needed. Without the violation of the explained requirement an extension of the developed stochastic solution approach considering more general precedence constraints between activities is straightforward. For this extension the variable~$z_{pqt}$ can again be the indicator variable for the decision that activity~$q$ of project~$p$ is scheduled in time period~$t$. A constraint $\sum_q z_{pqt} \le 1$ has to be added, and the precedence constraints have to be expressed by a more flexible set of constraints \citep{Artigues2008} as those given by Eq.~(\ref{eq:constraints3}).

\section{Solution procedures}
\label{sec:solutionprocedure}

\subsection{Problem structure}
\label{subsec:problemstructure}
As mentioned above, the considered project scheduling and staffing problem distinguishes between two decisions: the project scheduling decision where the execution times of the activities $q$ 
are determined, and the staffing decision where employees are assigned to cover parts of the work packages.
Let $z$ and $x$ denote the array of the decision variables $z_{pqt}$ for project scheduling, and the array of the decision variables $x_{ptsk}$ for staffing, respectively. By~$\xi$, we denote the random influence in our stochastic model. Abbreviating the total external cost by $G(z,x,\xi)$ (it depends on~$\xi$ since it is a random variable) and its expected value by $g(z,x)=E[G(z,x,\xi)]$, our optimization problem can be written as
$\min_{z,x} g(z,x)$,
where $(z,x)$ has to satisfy the constraints (\ref{eq:constraints2})~--~(\ref{eq:constraints6}).
For a given project schedule~$z$, we call
$\min_{x} \: g(z,x)$
on the constraint that $x$ is feasible with respect to (\ref{eq:constraints4})~--~(\ref{eq:constraints5}), the {\em staffing subproblem}.

\subsection{Expected value problem}
\label{subsec:meanvalueproblem}
The expected value (EV) problem or mean value problem is obtained from the original stochastic problem by replacing each random variable $D_{psq}$ by its expected value~$\bar{d}_{psq} = E(D_{psq})$, so that the distribution of $D_{psq}$ collapses to the point mass in $\bar{d}_{psq}$.
With the help of the introduction of the auxiliary variables for the external work time required $y_{pts}$ per project $p$ time period $t$ and skill $s$,
the following mixed integer linear problem is obtained:

\noindent\begin{tabularx}{\linewidth}{@{} L {80mm}L {20mm}>{\raggedleft\arraybackslash}X@{}}\\
\multicolumn{2}{c}{
\hspace*{-29ex} $\displaystyle \begin{aligned} \min_{z,x} \; \Theta_{EV}(z,x)
= \sum\limits_{p \in P} \: \sum\limits_{t = ES_p}^{LF_p} \: \sum\limits_{s \in S} \bigg ( c_s^e \, y_{pts}\bigg)\\
\end{aligned}$}
 & \tagarray\label{eq:EVobjectivefunction}\\
\\
subject to constraints \eqref{eq:constraints2}, \eqref{eq:constraints3}, \eqref{eq:constraints6}, and\\\addlinespace[2ex]
$\displaystyle \sum\limits_{q \in \uptau_{pt}}\bar{d}_{psq}z_{pqt} \leq y_{pts}+\sum\limits_{k \in K_s} x_{ptsk}$ &
$\displaystyle \begin{aligned}
	p &\in P,\\
	t&=ES_p,\dotsc,LF_p,\\
	s &\in S\\
\end{aligned}$
&\tagarray\label{eq:EVdemand}\\\addlinespace[2ex]

$\displaystyle y_{pts} \geq 0$ &
$\displaystyle \begin{aligned}
	p &\in P,\\
	t&=ES_p,\dotsc,LF_p,\\
	s &\in S\\
\end{aligned}$
&\tagarray\label{eq:constraints10}\\\addlinespace[2ex]
\end{tabularx}
The EV problem serves as an approximation for cases where the variances of the random variable $D_{psq}$ are very low, or, if the variances are larger, as a means to obtain an initial solution for an iterative search procedure


\subsection{Staffing subproblem}
\label{subsec:fixprojectschedule}
For a pre-defined feasible project schedule $z$,
the staffing problem $\min_{x} \: g(z,x)$ remains to be solved.
Because~$z$ is now already fixed,\\
\begin{tabularx}{\linewidth}{@{} L {80mm}L {20mm}>{\raggedleft\arraybackslash}X@{}}\\
\multicolumn{2}{c}{
\hspace*{18ex} $\displaystyle D'_{pst} = \sum\limits_{q \in \uptau_{pt}} D_{psq}z_{pqt}$
} & \tagarray\label{eq:D_pst}\\
\end{tabularx}
is a random variable that does not depend on any decision anymore, i.e., the distribution of~$D'_{pst}$ is already known.
$D'_{pst}$ represents the effective work time of project~$p$ in skill~$s$ that has been scheduled for time period~$t$ by the given schedule~$z$.
Using the random variables~$D'_{pst}$, we can express
the staffing problem in the form\\
\noindent\begin{tabularx}{\linewidth}{@{} L
{80mm}L {20mm}>{\raggedleft\arraybackslash}X@{}}\\
\multicolumn{2}{c}{
$\displaystyle\min_x \; \sum\limits_{p \in P} \: \sum\limits_{t = ES_p}^{LF_p} \: \sum\limits_{s \in S} c_s^e \: E \bigg( [D'_{pst}-\sum\limits_{k\in K_s} x_{ptsk}]^+\bigg)$
} & \tagarray\label{eq:objectivefunction2}\\
&
\end{tabularx}
subject to constraints \eqref{eq:constraints4} and \eqref{eq:constraints5}.


Obviously, the objective function \eqref{eq:objectivefunction2} is nonlinear, but it is not difficult to see that it is at least convex, since the function $x \mapsto x^+$ is a convex function and the expectation operator, as a linear operator, preserves convexity.
Let\\
\noindent\begin{tabularx}{\linewidth}{@{} L {80mm}L {20mm}>{\raggedleft\arraybackslash}X@{}}\\
\multicolumn{2}{c}{
$\displaystyle \varphi_{pst}(\zeta) =E([D'_{pst}-\zeta]^+)=\int_\zeta^\infty \!(\theta-\zeta) h_{pst}(\theta) \, \mathrm{d}\theta$
}&\tagarray\label{integral:integral15}\\\addlinespace[2ex]
&&
\end{tabularx}
with $h_{psq}(\theta)$ denoting the probability density function of $D'_{pst}$, such that (\ref{eq:objectivefunction2}) can be rewritten as\\
\noindent\begin{tabularx}{\linewidth}{@{} L {80mm}L {20mm}>{\raggedleft\arraybackslash}X@{}}\\
\multicolumn{2}{c}{
$\displaystyle \min_x \;
\sum\limits_{p \in P} \: \sum\limits_{t = ES_p}^{LF_p} \: \sum\limits_{s \in S}c_s^e \: \varphi_{pst} \Big (\sum\limits_{k\in K_s} x_{ptsk}\Big )$
}&\tagarray\label{integral:integral16}\\\addlinespace[2ex]
&&
\end{tabularx}
subject to \eqref{eq:constraints4} and \eqref{eq:constraints5}.
Elementary calculations show that\\
\noindent\begin{tabularx}{\linewidth}{@{} L {80mm}L {20mm}>{\raggedleft\arraybackslash}X@{}}\\
\multicolumn{2}{c}{
\hspace*{10ex}$\displaystyle \begin{aligned}
	\varphi'_{pst}(\zeta) &=\int_{-\infty}^\zeta \!h_{pst}(\theta)\mathrm{d}\theta-1=H_{pst}(\zeta)-1,\\
\end{aligned}$
}&\tagarray\label{integral:integral17}\\\addlinespace[2ex]
\end{tabularx}
with $H_{pst}$ denoting the cumulative distribution function (CDF) of $D'_{pst}$.

\subsection{Matheuristic}

Our first solution method is a {\em matheuristic} approach, i.e., a combination of a metaheuristic with an exact optimization technique. We solve the staffing subproblem by means of the exact {\em Frank-Wolfe} algorithm (see, e.g., \citet{Clarkson2010}) for convex optimization under linear constraints. This is the topic of Subsection~\ref{subsubsec:Frank Wolfe}. The scheduling problem is solved by a metaheuristic of Iterated Local Search type, which will be described in Subsection~\ref{subsubsec:metaheuristic}.

\label{subsec:matheuristic}
\subsubsection{Staffing: Frank-Wolfe algorithm}
\label{subsubsec:Frank Wolfe}
The idea of using the Frank-Wolfe algorithm for staffing problems has already been elaborated in~\cite{Gutjahr2013}.
Our technique used in the present work follows the approach described there rather closely, so we shall keep the presentation short.
First,
let us introduce the variables
$u_{ptsk}={x_{ptsk}}/{\eta_{sk}a_{kt}}$ $(p \in P, \: t =1,\ldots,T, \: s \in S, \: k \in K)$.
Using this substitution for a fixed project schedule $z$, the objective function~(\ref{eq:objectivefunction2}) can be rewritten as\\
\noindent\begin{tabularx}{\linewidth}{@{} L {80mm}L {20mm}>{\raggedleft\arraybackslash}X@{}}\\
\multicolumn{2}{c}{
$\displaystyle\min \; \Theta_z(u) \: = \: \sum\limits_{p \in P} \: \sum\limits_{t = ES_{p}}^{LF_{p}} \: \sum\limits_{s \in S}  c_s^e \: E \bigg( [D'_{pst}-\sum\limits_{k\in K_s} \eta_{sk}a_{kt}u_{ptsk}]^+\bigg)$.
} & \tagarray\label{eq:objectivefunction3}\\
&
\end{tabularx}

Moreover, we combine the project $p$ and the skill $s$ to the pair $\sigma=(p,s)$, which we call \textit{project-skill combination}, and we combine the time period $t$ and the resource $k$ to the pair $\nu=(t,k)$, the \textit{time-employee combination}.
Additionally, we re-label the pairs $\sigma=(p,s)$ and the pairs $\nu=(t,k)$ by introducing the new indices
$\sigma = 1,\dotsc,C = |P| \cdot |S|$ and
$\nu=1,\dotsc,L =|T| \cdot K$.
The variables~$u_{ptsk}$ are re-labelled accordingly, i.e., for $\sigma=(p,s)$ and $\mu=(t,k)$, the notation $u_{\sigma \nu}$ abbreviates $u_{ptsk}$.
In the constraints of Equation~\eqref{eq:constraints4}, the inequalities can be replaced by equalities, since there is an optimal solution in which each constraint~\eqref{eq:constraints4} is active.
Therefore, \eqref{eq:constraints4} can be replaced by
$\sum_{\sigma =1}^C u_{\sigma\nu} = 1 \quad (\nu=1,\ldots,L)$.
The column vector $u_{\nu}=(u_{1 \nu},\dotsc,u_{C \nu})'$ is an element of
the standard simplex in $\mathbb{R}^C$, hence the feasible set is a Cartesian product
of~$L$ standard simplices.

The idea of the iterative Frank-Wolfe algorithm is to replace in each iteration~$i$ the convex function $\Theta_z$ by its linear approximation at a current feasible solution~$u^{[i]}$. The approximating linear function has a minimizer $g^{[i]}$ on the feasible set which can be easily determined. Next, the minimum of the convex function $\Theta_z$ restricted to the line segment between $u^{[i]}$ and $g^{[i]}$ is identified. This can be done by line search. The minimizer found in this way is used as the new current solution for the next iteration.


For our application, the linear approximation in point~$u$ produces the optimization problem\\
\noindent\begin{tabularx}{\linewidth}{@{} L {80mm}L {20mm}>{\raggedleft\arraybackslash}X@{}}\\
\multicolumn{2}{c}{
$\displaystyle
\min\bigg \lbrace \sum\limits_{\sigma = 1}^C \sum\limits_{\nu = 1}^L \dfrac{\partial\Theta_z(u)}{\partial u_{\sigma\nu}}r_{\sigma\nu} \: \vert \: r= (r_{\sigma\nu})\in S_C^L \bigg \rbrace $.
}&\tagarray\label{eq:linearizationI}\\\addlinespace[2ex]
\end{tabularx}
This problem decomposes into $L$ partial problems\\
\noindent\begin{tabularx}{\linewidth}{@{} L {80mm}L {20mm}>{\raggedleft\arraybackslash}X@{}}\\
$\displaystyle
\min\bigg \lbrace \sum\limits_{\sigma = 1}^C \dfrac{\partial\Theta_z(u)}{\partial u_{\sigma\nu}}r_{\sigma\nu} \: \vert \: r_\nu= (r_{1\nu},\dotsc,r_{C\nu})' \in S_C\bigg \rbrace, $&
$\nu=1,\dotsc,L$.
&\tagarray\label{eq:linearizationII}\\\addlinespace[2ex]
\end{tabularx}

The solution of the $\nu$-th problem in~\eqref{eq:linearizationII} is an extremal point of the simplex $S_C$. Therefore, the point~$g^{[i]}$
is of the form $(e_{\sigma^*(1)},\dotsc,e_{\sigma^*(L)})$, where $e_\sigma$ is the $\sigma$-th unit vector, and $\sigma^*(\nu)$ is the index of the optimal
extremal point for time-employee combination~$\nu$. It is easily seen that $\sigma^*(\nu)$ is given by the index~$\sigma$ of the smallest value among the partial derivatives
in~(\ref{eq:linearizationII}).
One finds
$$
\dfrac{\partial \Theta_z(u)}{\partial u_{ptsk}}
= c_{s}^e \, \eta_{sk}a_{kt} \varphi'_{pst}\Bigg (\sum\limits_{k'\in K_s} \eta_{s k'} \, a_{k' t} \, u_{ptsk'}\Bigg)
$$
which, by (\ref{integral:integral17}), leads to the maximization problem

\noindent\begin{tabularx}{\linewidth}{@{} L {80mm}L {20mm}>{\raggedleft\arraybackslash}X@{}}\\
\multicolumn{2}{c}{
$\displaystyle
 \max_{p,s} \;  c_{s}^e \, \eta_{sk} \, a_{kt}\Bigg( 1-H_{pst}\bigg (\sum\limits_{k' \in K_s} \eta_{s k'}a_{k' t}u_{ptsk'}\bigg)\Bigg)$.
}&\tagarray\label{eq:minpartialderivate}\\\addlinespace[2ex]
&&
\end{tabularx}
Solving \eqref{eq:minpartialderivate} by enumeration, we find for each time-employee combination $\nu=(t,k)$ the project-skill combination $\sigma^*(\nu)=(p^*(\nu),s^*(\nu))$ providing the highest cost-reduction potential achievable by employee~$k$ in time period~$t$.

\begin{sloppypar}
Algorithm~\ref{alg:FrankWolfeAlgorithm1} presents the basic version of the Frank-Wolfe algorithm applied to the stochastic staffing subproblem.
We also slightly modify this basic procedure by doing the line search for each column~$\nu$ separately, followed by an immediate change of column~$\nu$ of the current matrix $u^{[i]}$. This gives Algorithm~\ref{alg:FrankWolfeAlgorithm2}.

\end{sloppypar}

\vspace*{-2ex}
\begin{algorithm}
\caption{Basic Frank-Wolfe algorithm for the staffing problem}
\label{alg:FrankWolfeAlgorithm1}
\begin{algorithmic}[1]
\State $i\gets 1$
\State $u^{[i]} \gets InitialSolution()$ \Comment{//Compute a feasible initial solution}
\Repeat
\State $i\gets i+1$
\ForAll {$\nu =(t,k)$}
\ForAll {$\sigma=(p,s)$ with $\sum_{q} z_{pqt} > 0$ and $s \in S_k$} 
\State $\chi(p,s) \gets c_s^e \, a_{kt} \, \eta_{sk} \, [1-H_{pst}(\sum\limits_{k' \in K_s} \eta_{s k'} \, a_{k' t} \, u^{[i]}_{ptsk'})]$ \label{alg:FWline:deriviation}
\EndFor
\State $\sigma^*(\nu) = (p^*(\nu), s^*(\nu)) \gets \arg\max_{p,s} \chi(p,s)$ \Comment{//Optimal index combination}
\State $g^{[i]}_\nu \gets e_{\sigma^*(\nu)}$\Comment{//Set the target column vector}
\EndFor
\State $g^{[i]} \gets (g^{[i]}_1,\ldots,g^{[i]}_L)$
\State $\vartheta^* \gets \arg\min_{0\leq\vartheta\leq1} \Theta_z ((1-\vartheta)u^{[i]} + \vartheta g^{[i]}$)\Comment{//Apply line search}
\State $u^{[i+1]} \gets (1-\vartheta^*)u^{[i]}+\vartheta^*g^{[i]}$ \Comment{//New $u^{[i]}$ according to $\vartheta^*$ }\label{alg:FWline:linesearch}
\Until{$i=i_{max}$}
\end{algorithmic}
\end{algorithm}

\vspace*{-6ex}
\begin{algorithm}
\caption{Modified Frank-Wolfe algorithm for the staffing problem}
\label{alg:FrankWolfeAlgorithm2}
\begin{algorithmic}[1]
\State $i\gets 1$
\State $u^{[i]} \gets InitialSolution()$ \Comment{//Compute a feasible initial solution}
\Repeat
\State $i\gets i+1$
\ForAll {$\nu =(t,k)$}
\ForAll {$\sigma=(p,s)$ with $\sum_{q} z_{pqt} > 0$ and $s \in S_k$} 
\State $\chi(p,s) \gets c_s^e \, a_{kt} \, \eta_{sk} \, [1-H_{pst}(\sum\limits_{k' \in K_s} \eta_{s k'} \, a_{k' t} \, u^{[i]}_{ptsk'})]$ \label{alg:FWline:deriviation}
\EndFor
\State $\sigma^*(\nu) = (p^*(\nu), s^*(\nu)) \gets \arg\max_{p,s} \chi(p,s)$ \Comment{//Optimal index combination}
\State $g^{[i]}_\nu \gets e_{\sigma^*(\nu)}$\Comment{//Set the target column vector}
\State obtain $g^{[i]}$ from $u^{[i]}$ by replacing column $u^{[i]}_\nu$ with $g^{[i]}_\nu$
\State $\vartheta^* \gets \arg\min_{0\leq\vartheta\leq1} \Theta_z ((1-\vartheta)u^{[i]} + \vartheta g^{[i]}$)\Comment{//Apply line search}
\State $u^{[i]} \gets (1-\vartheta^*)u^{[i]}+\vartheta^*g^{[i]}$ \Comment{//New $u^{[i]}$ according to $\vartheta^*$ }\label{alg:FWline:linesearch}
\EndFor
\State $u^{[i+1]} \gets u^{[i]}$
\Until{$i=i_{max}$}
\end{algorithmic}
\end{algorithm}
\vspace*{-3ex}


\subsubsection{Project Scheduling: Matheuristic solution method}
\label{subsubsec:metaheuristic}

\begin{figure}\centering
  \includegraphics[width=100mm]{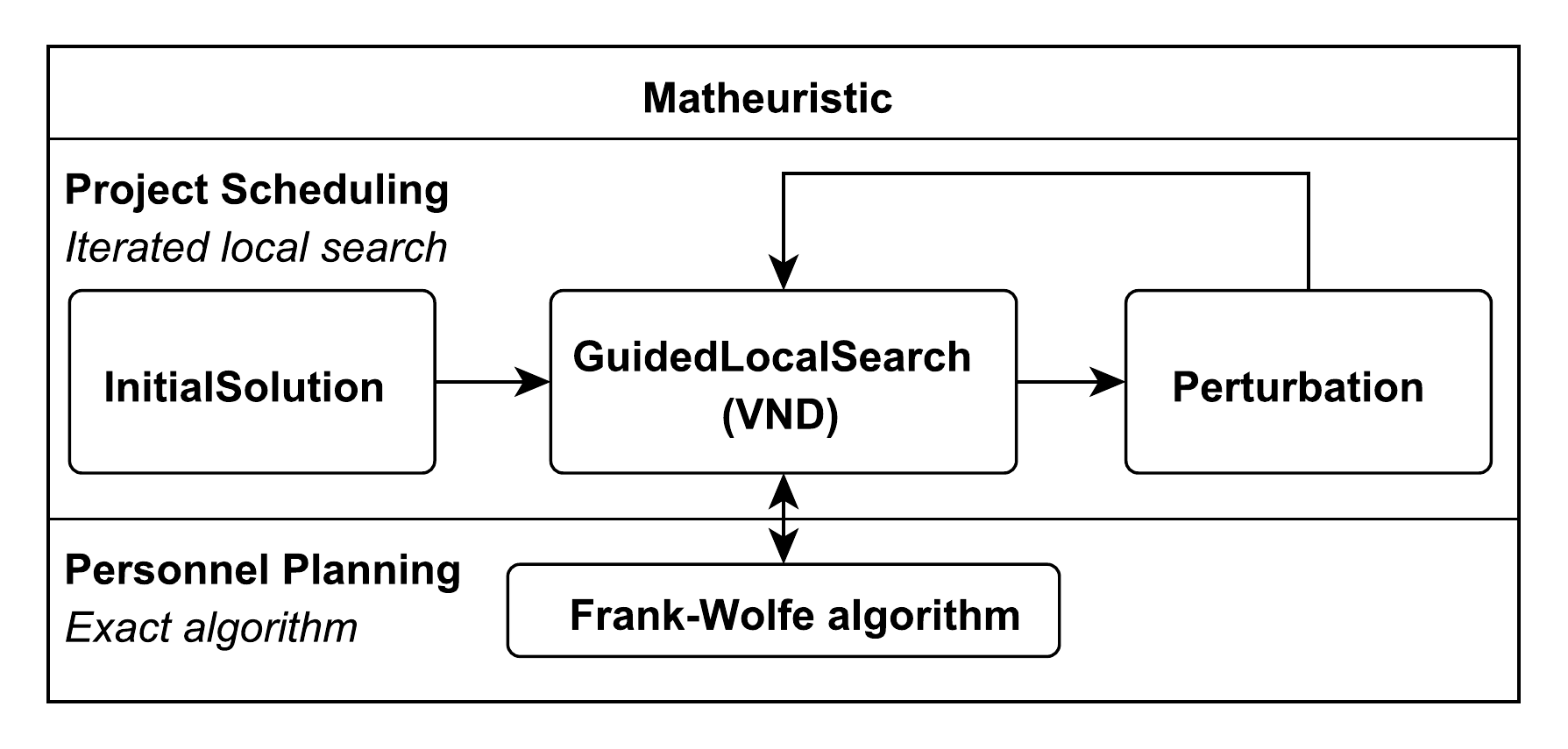}
\caption{Matheuristic: framework}
\label{fig:1}       
\end{figure}

After the description of the solution of the staffing subproblem in the previous section, the current section presents the search procedure that optimizes the project schedule~$z$. The implemented metaheuristic (which is also used for solving a deterministic version of our project scheduling problem in~\citet{Felberbauer2016b}) is a modification of iterated local search (ILS) \citep[e.g.,][]{Lourencco2010} using variable neighborhood descent (VND) as the local search component \citep[e.g.,][]{Hansen2010}.  The metaheuristic is combined with the exact Frank-Wolfe algorithm to a matheuristic (MH), schematically depicted in Figure~\ref{fig:1}.

\begin{sloppypar}
A pseudocode of the MH is given in Algorithm~\ref{alg:A2ILS}. Therein, instead of~$z$, an alternative representation of a project schedule by a two-dimensional scheme~$Z$ is used. The rows of $Z$ correspond to the projects $p$ $(1 \le p \le |P|)$, whereas the columns correspond to the activity indices $q=1,\ldots,d_p$ (since the numbers~$d_p$ need not to be identical, the number of entries per row can vary). The entry $Z_{pq}$ in row~$p$ and column~$q$ indicates the time period~$t$ in which activity~$q$ of project~$p$ is processed. A project schedule~$Z$ must satisfy the requirements on the time windows and on the starting time relationships, as defined in terms of the representation~$z$ by Equation~\eqref{eq:constraints2} and Equation~\eqref{eq:constraints3}, respectively. In Algorithm~\ref{alg:A2ILS}, some design parameters have already been set to fixed numerical values; these values resulted from numerous pre-tests as the most successful parametrizations.
\end{sloppypar}

\begin{algorithm}
\caption{MH search scheme}
\label{alg:A2ILS}
\begin{algorithmic}[1]
\State $n_r \gets 0$
\State $i_{min} \gets 50, \tilde{i} \gets 1000$
\State $Z, Z^* \gets InitialSolution()$
\While {$time<t_{max}$}
\State $k \gets 1$
\Repeat \Comment{//GLS-VND}
\State $n_r \gets n_r+1$
\State $i_{max} \gets n_r^{0.5} \cdot i_{min}$
\State $Z' \gets FirstImprovement(Z,Z^*,k,i_{max})$ \Comment{//Find the best neighbor}
\State $(Z^*,k) \gets neighborhoodChangeS(Z^*,Z',k,\tilde{i})$ \Comment{//Neighborhood change}
\State $Z \gets Z^*$
\Until{$k > k_{max}$}
\State $Z \gets Perturbation(Z, \: \lfloor \pi \cdot |P| \rfloor)$ \label{ILS:line:A2Perturbation}\Comment{//Perturbation}
\State $time \gets CpuTime()$
\EndWhile
\State\Return $Z^*$
\end{algorithmic}
\end{algorithm}


\textbf{Initial solution:} In the first step of MH, an initial solution~$Z$ for the project schedule
is generated by solving the Expected value problem presented in Section~\ref{subsec:meanvalueproblem}.
The solution~$Z$ is also used as the initialization of the incumbent solution~$Z^*$.

\begin{sloppypar}
\textbf{First improvement:}
The procedure \textit{FirstImprovement()} takes a current project schedule $Z$ and improves it by a local search to a solution $Z'$. The local search is based on a neighborhood defined by the current value of the parameter~$k$ (see below), and it is continued until (i) the neighbor solution is better than the current solution~$Z$,  (ii) the neighbor solution is better than the current incumbent~$Z^*$, or (iii) a local optimum is reached (no better neighbor solution exists). The name ``FirstImprovement'' indicates that we already terminate the search as soon as for the first time, an improving neighbor has been found; we do not necessarily explore the entire $k$-neighborhood.
The neighborhood definition, controlled by parameter~$k$, is based on a local move operator (``$k$-move'') that works as follows: First of all, for the given solution~$Z$, a capacity profile is determined. This capacity profile indicates the costs of external capacity needed for each time period.
Next, we identify that time period~$t$ for which the capacity costs are maximal, we consider the set of activities~$(p,q)$ that are executed in time period~$t$, and we sort these activities in descending order according to their contribution to the external capacity costs in period~$t$. A $k$-move to a neighbor solution consists in selecting $k$~consecutive activities~$(p,q)$ from the list and in shifting the execution time periods of the selected activities either one period forward or backward. Afterwards, in a subprocedure \textit{RepairSolution()}, all predecessor periods and successor periods of an affected project are checked for feasibility (see constraints in Equation~\eqref{eq:constraints3}) and, if necessary, repaired.
Finally, the objective function of the modified project schedule is evaluated, which requires a call of the Frank-Wolfe algorithm as explained in section~\ref{subsubsec:Frank Wolfe}. The numerical accuracy of the solution value determination by the Frank-Wolfe algorithm is controlled through a parameter~$i_{max}$, which will be explained in the Remark at the end of this subsection.
\end{sloppypar}

\textbf{Neighborhood change:} In the procedure~\textit{neighborhoodChangeS()}, the candidate solution $Z'$ that resulted from the local search in \textit{FirstImprovement()} is either accepted as the new incumbent~$Z^*$ or rejected, depending on its objective function value and on the outcome of a random event: If $Z'$ is better than the current~$Z^*$, then $Z'$ is accepted in any case, i.e., $Z^*$ is replaced by $Z'$, and the neighborhood size parameter~$k$ is reset to its initial value~1. Otherwise, we increase~$k$ by 1;  moreover, with probability~$\beta$, solution $Z'$ is accepted though it is worse than $Z^*$ (i.e., $Z^*$ is replaced again by $Z'$), and with probability $1-\beta$, the current incumbent solution~$Z^*$ is preserved. In the experiments, a value $\beta = 0.1$ turned out as a good choice.

The {\em repeat} loop in the algorithm performs a local search with a varying neighborhood definition, specified by the value of the parameter~$k$. Such a procedure is usually called Variable Neighborhood Descent (VND). For the maximal neighborhood size~$k_{max}$, we chose $k_{max} = 3$.


\textbf{Perturbation:} The procedure~$Perturbation()$
generates a new solution from the current solution~$Z$ by performing random swaps, using a parameter that indicates the number of swaps. The procedure is
triggered if all neighborhoods up to $k = k_{max}$ have failed to improve the incumbent solution. The number of swaps is chosen as equal to the product of a perturbation factor~$\pi$ and the number of projects $\vert P \vert$. Thus, we let the number of swaps depend on the problem size. A swap consists in exchanging the schedules (starting times of activities) of two randomly chosen projects and repairing the solution afterwards, if the new schedule is infeasible. In our experiments, it turned out that $\pi=0.15$ was a suitable value. Note that the perturbation procedure only works for projects with an identical number of activities.

{\bf Remark.}
When applying the Frank-Wolfe algorithm as a subprocedure of the heuristic algorithm~\ref{alg:A2ILS}, we adopt two strategies from~\citet{Gutjahr2013}.
First, in order to use computation time economically, the solution accuracy of the Frank-Wolfe algorithm can be increased gradually during the search process
as a function of the number~$n_r$ of conducted neighborhood searches. It turned out that using $i_{max}=i_{min} \cdot n_r^{0.5}$ iterations in the Frank-Wolfe algorithm, where $i_{min}$ is a fixed initial value, produced good results. For the evaluation in the procedure~\textit{neighborhoodChangeS()}, we always use a comparably large number~$\tilde{i}$ iterations.
In the second strategy, we suppose that the current upper bound $B(z,i)$ of the expected external costs after~$i$ iterations of the Frank-Wolfe algorithm converges with (approximately) exponential speed to the true value $a = \min_x g(z,x)$ as $i \to \infty$. That is, we assume $B(z,i) \sim a + b \cdot \exp(-ci)$, and estimate the parameters $a$, $b$, and $c$ from the three observations $B_1$, $B_2$ and $B_3$ of the bound in iteration $i=0.5 \cdot i_{max}$, $i=0.75 \cdot i_{max}$ and $i=i_{max}$, respectively.
This leads to the following estimate of the expected external costs:
$\tilde{\Theta}_z^*=a=(B_0B_2-B_1^2)/(B_0-2B_1+B_2)$.


\subsection{Sample average approximation}
\label{subsec:SAA}

\begin{sloppypar}
Our second solution approach for the stochastic project scheduling and staffing problem is the method of Sample Average Approximation (SAA), see \citet{kleywegt2002sample}. The SAA method samples scenarios from the given distribution of the random events and approximates the given stochastic problem by the so-called {\em deterministic equivalent}, the problem resulting as the average over the scenarios. The deterministic equivalent can be solved by methods from deterministic optimization.
\end{sloppypar}

In the case of our problem, we draw a set of $N$ random scenarios, described by the realizations $d_{psq}^{(1)}, \dotsc,d_{psq}^{(N)}$ of the random variables $D_{psq}$ $(p \in P, \: s \in S, \: 1 \le q \le d_p)$. Thus, $d_{psq}^{(n)}$ denotes the demand of work package~$(p,s,q)$ in scenario~$n$. The external work time required under scenario $n$ for project~$p$ and skill~$s$ in time period~$t$ will be denoted by the new variable $y_{pts}^{(n)}$. The deterministic equivalent is then the problem (\ref{eq:SAA})~--~(\ref{eq:SAAconstraint2}), which is obviously a mixed-integer linear program. Note that $\Theta_{SAA}(z,x)$ represents the average of the external costs over all scenarios, which approximates the expected external costs, and that the values of the variables $y_{pts}^{(n)}$ result in~(\ref{eq:SAAconstraint1}) by a formula analogous to~(\ref{eq:EVdemand}). The mixed-integer linear program (\ref{eq:SAA})~--~(\ref{eq:SAAconstraint2}) can be solved by standard solvers such as CPLEX. However, for realistic instance sizes and an appropriate number of samples, the problem (\ref{eq:SAA})~--~(\ref{eq:SAAconstraint2}) can become rather large, such that it may become difficult for the solver even to find a feasible solution within a reasonable time budget.


\noindent\begin{tabularx}{\linewidth}{@{} L {80mm}L {20mm}>{\raggedleft\arraybackslash}X@{}}\\
\multicolumn{2}{c}{
$\displaystyle \begin{aligned}
\hspace*{-22ex} \min_{z,x} \;
\Theta_{SAA}(z, x) = \frac{1}{N} \: \sum\limits_{n \in N} \: \sum\limits_{p \in P} \: \sum\limits_{t = ES_p}^{LF_p} \: \sum\limits_{s \in S} c_s^e \, y_{pts}^{(n)}&\\
\end{aligned}$
} & \tagarray\label{eq:SAA}\\
\\
subject to constraints (\ref{eq:constraints2})~--~(\ref{eq:constraints6}) and\\\addlinespace[2ex]
$\displaystyle \sum\limits_{q \in \uptau_{pt}}d_{psq}^{(n)}z_{pqt} \leq y_{pts}^{(n)}+\sum\limits_{k \in K_s} x_{ptsk}$ &
$\displaystyle \begin{aligned}
	p &\in P,\\
	t&=ES_p,\dotsc,LF_p,\\
	s &\in S,\\
	n &\in N\\
\end{aligned}$
&\tagarray\label{eq:SAAconstraint1}\\\addlinespace[2ex]

$\displaystyle y_{pts}^{(n)} \geq 0$ &
$\displaystyle \begin{aligned}
	p &\in P,\\
	t&=ES_p,\dotsc,LF_p,\\
	s &\in S,\\
	n &\in N\\
\end{aligned}$
&\tagarray\label{eq:SAAconstraint2}\\\addlinespace[2ex]
\end{tabularx}

\section{Experimental results}
\label{sec:experiments}

\subsection{Test instance generation}
\label{subsec:Testinstances}

For our computational experiments, we generated test instances by the test instance generator proposed in~\cite{Heimerl2010a}. This test instance generator is inspired by data from the IT department of a large semiconductor manufacturer. However, we had to extend the test instances since in the present paper, the information about the demand of work packages is assumed to be uncertain. In particular, in addition to the parameters number of projects~$|P|$, time window size~$\gamma=LS_p-ES_p$, and number of skills per resource~$\vert S_k \vert$, also a further parameter representing the degree of uncertainty (which will be explained later) had to be varied across the test instances. The database with the instance set is available for download at the website  \url{http://phaidra.fhstp.ac.at/o:2529}.

Table \ref{tab:basescenario} lists the parameter values for the basic instance structure. For this group of instances, the number of projects is~10, and the earliest start period~$ES_p$ of each project $p$ is drawn from a uniform distribution between 1 and 7. The project length~$d_p$ is set to six periods for each project, and the planning horizon is defined as $T=12$, which represents annual strategic project scheduling and personnel planning. We choose a time window size $\gamma=1$ which means that for each activity $q$, there are two possible start times available, considering precedence constraints. With $S^{(p)}$ denoting the set of skills required by project~$p$, and with $S^{(p,q)}$ denoting the set of skills required by activity~$q$ of project~$p$, the test instance generator limits the number $|S^{(p)}|$ of skills per project by a pre-defined bound; moreover, for each activity~$q$ of a project~$p$, the test instance generator specifies the number~$|S^{(p,q)}|$ of required skills. For all skills~$s$ not contained in $S^{(p,q)}$, the values $d_{psq}$ are set to zero. (For details concerning this aspect of the test instance generation, see~\cite{Heimerl2010a}). In our case, we chose $|S^{(p,q)}|= 2$ and limited the total number of skills per project by $|S^{(p)}|\leq 3$. Ten internal resources are assumed ($|K|=10$), each owning $|S_k|=2$ out of $|S|=10$ skills. The resources have different efficiency values~$\eta_{sk}$ for each skill~$k$ they own; these efficiency values are drawn from a truncated normal distribution with an expected value of $\mu=1$, a standard deviation of $\sigma=0.25$, and minimum and maximum threshold values of 0.5 and 1.5, respectively. The available capacity per internal resource $a_{kt}$ is 20 per time period. Note that we do not assume that a time period has a length of one time unit: In our above-mentioned interpretation of the test instances as referring to annual planning, the time unit is a day, and a time period extends over a month, so $a_{kt}=20$ means that an employee works 20 days per month.

The external cost rates $c_s^e$ differ for different skills and are drawn from a truncated normal distribution $TN_{a,b}(\mu,\sigma)$ with $\mu=800$, $\sigma=100$, $a=600$ and $b=1000$. In test instance generation, (planned) utilization~$\rho$ is defined as the ratio of the overall expected resource demand to available internal resource capacity. In the basic instance structure, we set~$\rho=1$. For each work package, the test instance generator computes an initial value $[E(D_{psq})]^{init}$ of the expected resource demand from the utilization~$\rho$ and the resource supply values~$a_{kt}$. The actual expected value of the resource demand $E(D_{psq})$ is then drawn from a normal distribution, with mean $\mu=[E(D_{psq})]^{init}$ and a coefficient of variation $CV =0.1$. For the basic instance structure, we assume a symmetric triangular distribution of the {\em actual} demand with parameters $(D_{psq}^{min},D_{psq}^{mod},D_{psq}^{max})$, where $D_{psq}^{mod}=E(D_{psq})$, $D_{psq}^{max}=D_{psq}^{mod} \cdot c_{max}$, and $D_{psq}^{min}=D_{psq}^{mod} \cdot c_{min}$. The interval $[c_{min},c_{max}]$ controls the level of uncertainty; if $c_{min} = c_{max}$, we get the deterministic boundary case. For the basic instance structure, a moderate level of uncertainty is assumed by setting $c_{min}=0.7$ and $c_{max}=1.3$. Notice that because we use {\em symmetric} triangular distributions in our basic test instances, the expected value $E(D_{psq})$ is identical to the modal value $D_{psq}^{mod}$ of the distribution.\footnote{The choice of the triangular distribution for related models is discussed, e.g., in \citep{Law1991,Wing1995}. \citet{Law1991} see the triangular distribution as a rough model suitable for cases where there is limited data available or the costs of data collection are high. Also the Beta distribution would be a candidate for modeling the distributions of the random variables $D_{psq}$. However, as shown in~\citet{Gutjahr2013} for a related model, using a Beta distribution within the Frank-Wolfe framework comes at the price of a distinct increase in computation time, and the results are typically similar to those obtained by the triangular distribution.}

\begin{table}\centering
\caption{Basic instance structure}
\label{tab:basescenario}
\begin{tabular}{ll}
\hline
$ \vert P\vert =10$ & $|S^{(p)}| \leq 3$\\
\hline
$ ES_p \sim U(1,7)$ & $|K|=10$ \\
\hline
$ d_p=6$ & $\vert S_k \vert=2$ \\
\hline
 $T=12$ & $\eta_{sk} \sim TN_{0.5,1.5}(1,0.25)$ \\
\hline
 $\gamma=LS_p-ES_p=1$ & $a_{kt}=20$ \\
\hline
 $|S| = 10$ & $c_s^e \sim TN_{600,1000}(800,100)$\\
\hline
$|S^{(p,q)}|=2$& $\rho = 1.0$\\
\hline
$[c_{min},c_{max}]=[0.7, 1.3]$\\
\hline
\end{tabular}
\end{table}

The basic instance structure is varied then to obtain other instance structures, according to Table~\ref{tab:experimentallevels} which lists the parameters with their used values.
For test instance generation, we use a {\em ceteris paribus design}, which means that we fix all parameters on the value of the basic instance structure and vary the value of one investigated parameter. This yields $4+(4-1)+(6-1)+(3-1)=14$ different instance structures.
For each instance structure, we generate 10 instances,
which leads to 140 test instances.  Additionally, for each of these 150 test instances, four different levels of the degree of uncertainty (see Table~\ref{tab:aou}), are investigated. This produces a total of 560 instances.

All tests are performed on a standard PC with an Intel Quad Core Processor. In detail, we used an Intel Xeon E3-1271 v3 processor (Frequency: 3,60GHz) with eight kernels and 32 Gigabytes working memory. All presented algorithms are implemented in Eclipse using Java version 1.7. We use the Java API of ILOG CPLEX version 12.4 for the SAA and the EV model formulations.

\begin{table}\centering
\caption{Parameters and experimental values for the test instance generation}
\label{tab:experimentallevels}
\begin{tabular}{ll}
\hline
\textbf{Parameter} & \textbf{Experimental value} \\
\hline
$ \vert P \vert$ & $10,15,20,25$\\
\hline
$\gamma$=$LS_p-ES_p$ & $0,1,2,3$\\
\hline
$\vert S_k \vert$ & $1,2,4,6,8,10$\\
\hline
$\rho$ & $0.8,1.0,1.2$\\
\hline
\end{tabular}
\end{table}

\begin{table}\centering
\caption{Investigated degree of uncertainty}
\label{tab:aou}
\begin{tabular}{ll}
\hline
\textbf{Parameter} & \textbf{Experimental values} \\
\hline
$[c_{min},c_{max}]$ & $[0.9,1.1], \: [0.7,1.3], \: [0.5,1.5], \: [0.2,1.8]$\\
\hline
\end{tabular}
\end{table}

\subsection{Parameter setting for the Frank-Wolfe algorithm}
\label{subsec:Frank Wolfe setting}

There are several design decisions that have to be made when using the Frank-Wolfe algorithm within our matheuristic framework procedure: First of all, two slightly different implementation variants of this algorithm were presented in Section~\ref{subsubsec:Frank Wolfe}. Secondly, the algorithm requires an initial solution; we consider two options for its choice. Finally, it has to be specified how the line search is done; again, two different options will be investigated.

\textbf{\textit{Algorithmic variants:}} We shall compare the basic Algorithm~\ref{alg:FrankWolfeAlgorithm1} (``old'' in Table~\ref{tab:FrankWolfeSetting}) to the modified Algorithm~\ref{alg:FrankWolfeAlgorithm2} (``new'' in Table~\ref{tab:FrankWolfeSetting}).

\textbf{\textit{Initial solution:}} For the determination of the initial solution, two alternative approaches are tested. First, we use the MILP-solver CPLEX (``lp'' in Table~\ref{tab:FrankWolfeSetting}) to solve the deterministic linear staffing problem defined by the parameters $d'_{pst}=E(D'_{pst}) = \sum_q \bar{d}_{psq} \, z_{pqt}$. Secondly, we apply a greedy staffing heuristic~\citep{Felberbauer2016b} to solve the same deterministic counterpart problem heuristically (``gh'' in Table~\ref{tab:FrankWolfeSetting}). 

\textbf{\textit{Line search method:}}
For the line search step of the algorithm, two methods are analyzed: The first method (``gs'' in Table~\ref{tab:FrankWolfeSetting}) uses Golden Section Search according to~\citet{Kiefer1953}.
The second method (``fs'' in Table~\ref{tab:FrankWolfeSetting}) follows the
suggestion in~\citet{Clarkson2010}: it refrains from determining the $\arg\min$ in line~13 of Algorithm~\ref{alg:FrankWolfeAlgorithm1} or line~12 of Algorithm~\ref{alg:FrankWolfeAlgorithm2}, respectively, but uses instead in each iteration~$i$ a pre-defined step size $\vartheta^* = \vartheta^*(i)$ depending on the iteration index. The value of $\vartheta^*$ is calculated as
$\vartheta^*=2/(i+2)$ $(i=1,2,\ldots)$.
It is clear that the value of $\vartheta^*$ determined in this way does not produce the minimizer on the line segment between~$u^{[i]}$ and $g^{[i]}$, but by the special choice of the step sizes (convergence to zero and finiteness of the sum of the squares), the convergence property of the Frank-Wolfe algorithm
to the exact overall minimizer of $\Theta_z$
is preserved (for details, see \citet{Clarkson2010}). The advantage of the fixed step sizes scheme is that it does not require an evaluation of the function values $\Theta_z$ during the execution of the algorithm; it suffices to evaluate the derivatives of $\Theta_z$.



Combining the two alternative options for each of the three design decisions indicated above, we get $2^3=8$ different design variants of the Frank-Wolfe algorithm. The following results compare the performance of these eight design variants. 
For each design variant, we shall report
its average solution value~$sv$ and its average computation time~$ct$ at $1000$ randomly selected time schedules $z=(z_{pqt})$ for a single fixed problem instance generated according to the basic instance structure.
In a pre-test, it turned out that a number~$i_{max}=1000$ of iterations was sufficient to get close enough to the value $\Theta_z$ achieved by a much higher number~$10^6$ of iterations.
Therefore, we used $i_{max}=1000$ for all design variants.


\begin{table}
\caption{Average solution value \textit{(sv)}, standard deviation of the solution value, average computation time \textit{(ct)}, and standard deviation of the computation time for the eight Frank-Wolfe design variants. Last two columns: relative values compared to those of the best design variants. Based on  $i_{max}=1000$ iterations and $1000$ schedules $z$ of a fixed instance.}
\label{tab:FrankWolfeSetting}
\begin{tabular}{p{2.5cm}ccrr p{0.01cm}cr}
\hline
\multirow{3}{*}{\parbox{1.5cm}{\textbf{FW design variant}}} & \multicolumn{4}{c}{\textbf{absolute}}&& \multicolumn{2}{c}{\textbf{relative}} \\
&\textit{sv}&\textit{sv-stdev}&\textit{ct}&\textit{ct-stdev}&&$\Delta sv$&$\Delta ct$\\
&[CU]&[CU]&[ms]&[ms]&&[\%]&[\%]\\
\hline
\cmidrule{2-5}\cmidrule{7-8}
\textit{1. lp-fs-new} & $286,341.3$&$0.00$ &$102$&$6$ &&$0.00$& $3.15$\\
\textit{2. lp-fs-old} & $286,341.4$&$0.00$ &$213$&$14$ &&$0.00$& $115.00$\\
\textit{3. lp-gs-new} & $286,342.1$&$0.00$ &$22,439$&$118$ &&$0.00$& $22,515.00$\\
\textit{4. lp-gs-old} & $286,354.5$&$0.00$ &$14,448$&$112$ &&$0.00$& $14,462.00$\\
\textbf{\textit{5. gh-fs-new}} & $\textbf{286,341.2}$&$\textbf{0.00}$ &$\textbf{99}$&$\textbf{53}$ &&$\textbf{0.00}$& $\textbf{0.00}$\\
\textit{6. gh-fs-old} & $286,341.3$&$0.00$ &$209$&$9$ &&$0.00$& $110.00$\\
\textit{7. gh-gs-new} & $286,498.4$&$0.00$ &$30,015$&$200$ &&$0.05$& $30,151.00$\\
\textit{8. gh-gs-old} & $286,554,8$&$0.00$ &$15,095$&$88$ &&$0.07$& $15,113.00$\\
\hline
\end{tabular}
\end{table}

\begin{sloppypar}
The results are shown in Table~\ref{tab:FrankWolfeSetting}. The comparison between the two algorithmic variants Algorithm~\ref{alg:FrankWolfeAlgorithm1} and the new Algorithm~\ref{alg:FrankWolfeAlgorithm2} (immediately applying the best partial derivative for the update of the project plan) shows a slight superiority of Algorithm~\ref{alg:FrankWolfeAlgorithm2}. For the decision on the used initial solution method, the test shows that the computation time for solving the deterministic staffing problem takes in the average $\approx 3.73$ ms using the LP-solver and $\approx 0.26$ ms using the greedy heuristic. On the other hand, the expected cost of the initial staffing plan $x$ in iteration $i=1$ according to the greedy heuristic is $\approx9\%$ higher than the one obtained by the exact LP solver. Nevertheless, applying a two-tailed sign test to the final results for solution values and computation times shows no significant difference between the performance of the LP solver and the greedy heuristic (significance level $\alpha = 0.05$). For larger instances, where the LP solving time increases rapidly, the greedy heuristic can become the only feasible alternative, so that in total, we may give a preference to the greedy heuristic. Concerning the line search method, finally, it can be seen that the pre-defined step sizes scheme clearly outperforms the Golden Section Search: Although Golden Section Search provides faster improvements in the first few iterations than the step sizes scheme, the solution quality after 1000 iteration is not better, and the computation time per iteration is $\approx 200$ times higher.

\end{sloppypar}

Summarizing, the best-performing Frank-Wolfe variant uses the immediate application of each best partial derivative for the update of the project plan, the greedy heuristic for the calculation of the initial staffing solution, and the pre-defined step sizes scheme.

\subsection{Parameter setting for Sample Average Approximation}
\label{subsec:SAAPara}

The crucial parameter of the Sample Average Approximation procedure is the number~$N$ of sampled random scenarios. Therefore, the following pre-tests have been conducted to find an appropriate value of~$N$. It is clear that the objective function of~(\ref{eq:SAA}) is only an approximation to the true objective function, such that even if the SAA problem is solved exactly, we do not necessarily obtain the exact solution of the original problem. To explore the tradeoff between the two effects of increasing the value of~$N$, namely to improve the accuracy of the objective function estimation on the one hand, and to increase the computation time on the other hand,
a subset of our instances has been investigated. We chose the instances of those instances structures where the number of projects is varied as $|P| \in \lbrace 10,15,20,25 \rbrace$, the time window size is $\gamma \in \lbrace 0,1,2,3 \rbrace$,  and the number of skills per resource as well as the utilization are fixed to the values of the basic instance structure, i.e., $|S_k|=2$ and $\rho=1$. The degree of uncertainty was varied as $[c_{min}, \, c_{max}] \in \lbrace [0.9,1.1],\: [0.7,1.3],\: [0.5,1.5],\: [0.2,1.8] \rbrace$. For these instances, we varied the sample size as~$N \in \lbrace 10,20,30,\dotsc,100 \rbrace$ and analyzed the solution time and the achieved solution quality.
By the SAA model from section~\ref{subsec:SAA} with sample size~$N$, we compute the solution $(z_{SAA}^*,x_{SAA}^*)$, where $z_{SAA}^*$ is the optimal project plan and $x_{SAA}^*$ is the optimal staffing plan. Now, we compare the obtained objective function value $\Theta_{SAA}(z_{SAA}^*,x_{SAA}^*)$ of~(\ref{eq:SAA}) to the {\em true} evaluation $E[G(z_{SAA}^*,x_{SAA}^*,\xi)]$ of the solution $(z_{SAA}^*,x_{SAA}^*)$ according to the underlying exact probability model (cf.~the notation in section~\ref{subsec:problemstructure}). The relative gap between the two evaluations is described by

\noindent\begin{tabularx}{\linewidth}{@{} L {80mm}L {20mm}>{\raggedleft\arraybackslash}X@{}}\\
\multicolumn{2}{c}{
$\displaystyle|E[G(z_{SAA}^*,x_{SAA}^*,\xi)] - \Theta_{SAA}(z_{SAA}^*,x_{SAA}^*)|/E[G(z_{SAA}^*,x_{SAA}^*,\xi)]$.
}&\tagarray\label{eq:SAAvsESAA}\\
&
\end{tabularx}



\begin{figure}
\hspace*{-7ex}
{\includegraphics[scale=0.46]{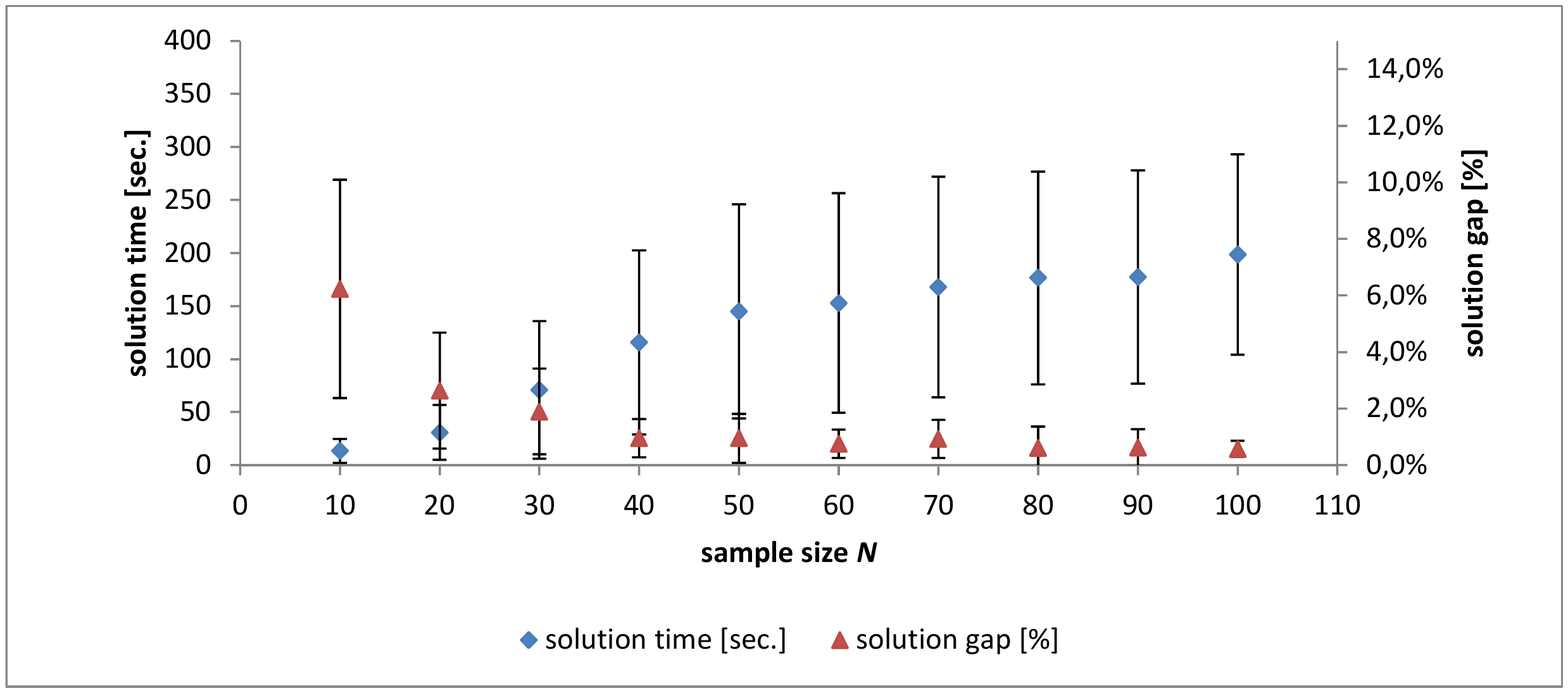}}
\vspace*{-22ex}
\caption{Performance of solution time and solution gap in dependence of sample size $N$ for
$\vert P \vert=20$, $\gamma=1$, $\rho=1.0$, $ \vert S_k\vert=2$, and $[c_{min}, \, c_{max}]=[0.7,1.3]$.}
\label{fig:SampleSize}
\end{figure}

In Figure~\ref{fig:SampleSize}, the solution time and the relative solution gap according to Eq.~(\ref{eq:SAAvsESAA}) as well as their 95\% confidence intervals are depicted in dependence of the sample size $N$. We show here the special case of $\vert P \vert=20$ projects and the other parameters as in the basic instance structure.
It can be observed that
for a sample size of $N=100$, the solution gap is $\approx 0.5\%$, i.e.,
the average objective function value over the scenarios can be considered as a good estimate for the expected external costs.
A further observation is that the solution time of the SAA model varies to a considerable extent.
This behavior points out a first drawback of relying on the exact solution of the SAA model to solve our problem.

These results were extended by solving all test instances of the instance subset specified at the beginning of this subsection. With an appropriate sample size of $N=100$ and a computation time limit of 360 sec for the CPLEX solver, we observed that CPLEX was able to solve the SAA problems for all considered test instances. In the average over all considered test instances, a relative gap of $0.71\%$ according to Eq.~(\ref{eq:SAAvsESAA}) was obtained.

\subsection{Matheuristic vs. Sample Average Approximation}
\label{subsec:SAAvsILS}

\begin{figure}
\hspace*{-5ex}
\includegraphics[trim=25mm 50mm 25mm 42mm,clip,scale=0.48]{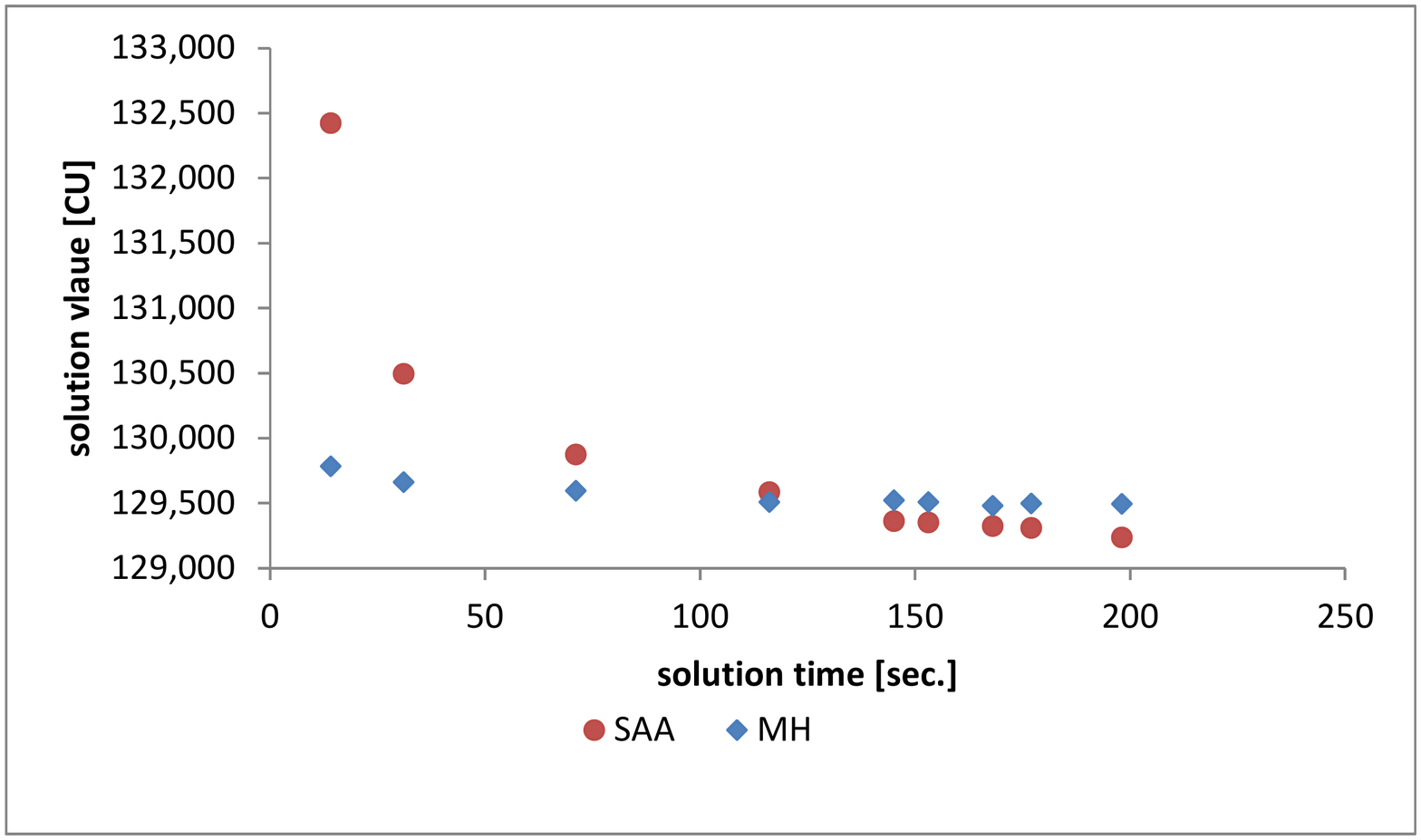}
\caption{Solution value performance of the sample average approximation approach (SAA, circles) and the matheuristic (MH,  diamonds) with respect to the given time budget for the instance structure with $\vert P \vert =20$, $\gamma=1$, $\rho=1.0$, $ \vert S_k\vert=2  $, and $[c_{min}, \: c_{max}] = [0.7,1.3]$.}
\label{fig:SAAvsILS}
\end{figure}

This section reports on the numerical comparison of the two presented solution methods, i.e., the developed matheuristic (MH) and the Sample Average Approximation (SAA) approach.
Note that by fixing a sample size $N$, the computation time consumed by the SAA approach
is already defined.
To ensure a fair comparison between SAA and MH, we computed, for a given test instance structure, the average computation times of the SAA approach for sample sizes $N = 10, 20,\ldots, 100$
and used each of these ten values as the time budget (termination criterion) for a corresponding run of the MH approach.
For each problem instance structure, ten random instances were generated, and for each of these generated instances, ten optimization runs (with different seed values for the random number generator) were executed. This led to 100 solution values for MH and 100 solution values for SAA per time budget. The averages of these solution values for a special instance structure are depicted in Figure~\ref{fig:SAAvsILS}. The solution values have been computed based on the determination of the exact objective function value of the solution $(z,x)$ provided by the respective approach, i.e., the value $E[G(z,x,\xi)]$. The reader will observe that Figure~\ref{fig:SAAvsILS} only contains 9 solution time values instead of 10, as one would expect. This is because sample sizes $N=80$ and $N=90$ led to identical solution times, cf.~Figure~\ref{fig:SampleSize}.

The results show that the solution quality of MH is less sensitive to the time budget than that of SAA.
Applying a two-tailed sign test to the results of each of the ten different time budgets, we found that the MH results were significantly superior, at significance level $\alpha=0.05$, for the smallest time budget $14$ $sec$ which was the time budget produced by sample size $N=10$. For larger time budgets, the sign test could not confirm statistically significant superiority at level $\alpha=0.05$ of either SAA over MH or vice versa. However, this lack of significance may be due to our small sample size of ten instances (in order to ensure independence, we had to take the average over the ten runs of each instance).
Summarizing, in the considered instances structure, MH and SAA produce results of a comparable quality in the medium computation time range, with a slight but nonsignificant advantage for SAA. For small computation times, MH is superior.

Next, we compare SAA, applying sample size $N=100$, to MH, using for both the same time limit of $360$ seconds. Let $(z_{MH}^*,x_{MH}^*)$ and $(z_{SAA}^*,x_{SAA}^*)$ denote the optimal project schedule and staffing plan according to MH and SAA, respectively. To do the comparison, we compute the normalized difference\\
\noindent\begin{tabularx}{\linewidth}{@{} L {80mm}L {20mm}>{\raggedleft\arraybackslash}X@{}}\\
\multicolumn{2}{c}{
 $(E[G(z_{SAA}^*,x_{SAA}^*,\xi)] - E[G(z_{MH}^*,x_{MH}^*,\xi)])/E[G(z_{SAA}^*,x_{SAA}^*,\xi)]$.
 }&\tagarray\label{eq:SAAvsMH}\\
&
\end{tabularx}
of the expected external costs. Averaging this measure over all test instances, we obtained a value of $0.1288$, indicating that in the average, the SAA solutions are by $12.88\%$ worse than the MH solutions. However, this result should be interpreted cautiously. The considerable difference can mainly be attributed to test instances where the expected external costs of the optimal solutions are almost zero. 
Such a situation easily occurs in instances where the number $\vert S_k \vert$ of skills per resource is high. Investigating the same measure only for the basic instance structure, we find that there is no significant difference between the two solution methods.

Real-world project scheduling problems are often rather large. To investigate which effect an increasing size of the problem instance exerts on the comparison between MH and SAA, we generated test instances where both the number of projects and the number of resources were increased simultaneously as $|P|=|K|=50,100,150,200,250$, and the parameter values $\gamma=2$, $|S_k|=2$ and $\rho=1$ from Table~\ref{tab:basescenario} were used.
For the SAA model, we applied CPLEX with default values and set the time limit to 10 hours, which is an appropriate time budget for a tactical planning problem. We found that the solver could not provide a feasible integer solution for $40\%$ of these instances; the share of solvable instances rapidly drops as $|P|$ and $|K|$ become larger than $150$. However, even for the instances for which the SAA model can be feasibly solved, the SAA solution is in the average by $45.57\%$ worse than the MH solution. Additionally, MH needs only a fraction of the SAA solution time to find a good solution.

We conclude that the developed matheuristic is a robust solution procedure that performs well
for small and medium-sized test instances, and offers good solutions also where the SAA formulation fails to return a feasible solution. The main drawback of the SAA approach is its poor reliability with respect to the identification of a feasible solution and its volatile solution time. Nevertheless, 
for not too large instances, the SAA approach is a promising alternative to the application of a (partially) heuristics-based method, especially in cases where powerful hardware is available.

\subsection{Deterministic vs. stochastic planning}
Is the advantage of treating the given project scheduling and staffing method as a stochastic optimization problem substantial enough to justify
the increased computational effort, compared to a simplified, deterministic formulation? To shed light on this question, we deepen our experimental analysis in the following two subsections.
In Section~\ref{subsubsec:accuracyEV}, the accuracy of using the solution value of the expected value problem as a forecast for budget planning is investigated, whereas in Section~\ref{subsubsec:VSS}, the value of the stochastic solution is discussed.

\subsubsection{Accuracy of the expected value problem}
\label{subsubsec:accuracyEV}

\begin{sloppypar}
The solution of the expected value (EV) problem (\ref{eq:EVobjectivefunction})~--~(\ref{eq:constraints10}) provides a manager with a project schedule $z_{EV}^*$, a staffing plan $x_{EV}^*$, and a corresponding solution value $\Theta_{EV}(z_{EV}^*,x_{EV}^*)$.
The obtained forecast $\Theta_{EV}(z_{EV}^*,x_{EV}^*)$ of the external costs could be used as an input for the budget planning process.
If the demand information is actually stochastic, this forecast will be rather rough, and it will tend to underestimate the true costs.
Evaluating the obtained project and staffing plan based on the stochastic model by computing $E[G(z_{EV}^*,x_{EV}^*,\xi)]$ gives a hint of how the forecast provided by solution of the EV problem will perform in the stochastic environment. The difference\\
\noindent\begin{tabularx}{\linewidth}{@{} L {80mm}L {20mm}>{\raggedleft\arraybackslash}X@{}}\\
\multicolumn{2}{c}{
$\displaystyle BB=E[G(z_{EV}^*,x_{EV}^*,\xi)]-\Theta_{EV}(z_{EV}^*,x_{EV}^*)$
}&\tagarray\label{eq:SAAvsESAAabs}\\
&
\end{tabularx}
could be considered as a {\em budget bias} caused by a deterministic solution approach. Table \ref{tab:relativeaccurarcy} shows the relative bias $BB^{rel}$ of the EV problem solution value, i.e., the quotient\\
\noindent\begin{tabularx}{\linewidth}{@{} L {80mm}L {20mm}>{\raggedleft\arraybackslash}X@{}}\\
\multicolumn{2}{c}{
$\displaystyle	BB^{rel}=(E[G(z_{EV}^*,x_{EV}^*,\xi)]-\Theta_{EV}(z_{EV}^*,x_{EV}^*))/E[G(z_{EV}^*,x_{EV}^*),\xi]$,
}&\tagarray\label{eq:SAAvsESAArel}\\
&
\end{tabularx}
for the basic instance structure and for all instances, in dependence of the degree of uncertainty.
\begin{table}\centering
\caption{Relative accuracy $BB^{rel}$ of the EV problem, measured in \%, (i) for the basic instance structure and (ii) for all instances, in dependence of the degree of uncertainty}
\label{tab:relativeaccurarcy}
\begin{tabular}{llllll}
\hline
\multirow{2}{*}{\parbox{1.5cm}{Instances}} & \multicolumn{5}{c}{$[c_{min},c_{max}]$}\\
&[0.9,1.1]& [0.7,1.3]& [0.5,1.5]& [0.2,1.8]& $\varnothing$ \\
\hline
Basic instance structure&3.96\%&11.32\%&17.94\%&26.68\%&14.97\%\\
\hline
Avg. over all instances&33.33\%&41.76\%&47.58\%&54.20\%&44.22\%\\
\hline
\end{tabular}
\end{table}
We see that with increasing uncertainty, the bias increases, and that it reaches fairly large values. For the basic instance structure, already under a moderate level of uncertainty of $[c_{min},c_{max}] = [0.7,1.3]$, the use of a deterministic planning approach leads to an underestimation of the external costs by 11\%. Averaged over the instances from all instance structures (with ten random instances from each instance structure), this deviation is even distinctly higher (42\%). According to these results, we conclude that the deterministic EV approach to the considered project scheduling and staffing problem leads to a systematic underestimation of external costs and can, as a consequence, seriously threaten the budget plan.
\end{sloppypar}

\subsubsection{Value of the stochastic solution}
\label{subsubsec:VSS}

\begin{figure}\centering
  \includegraphics[trim=25mm 42mm 25mm 42mm,clip,scale=0.35]{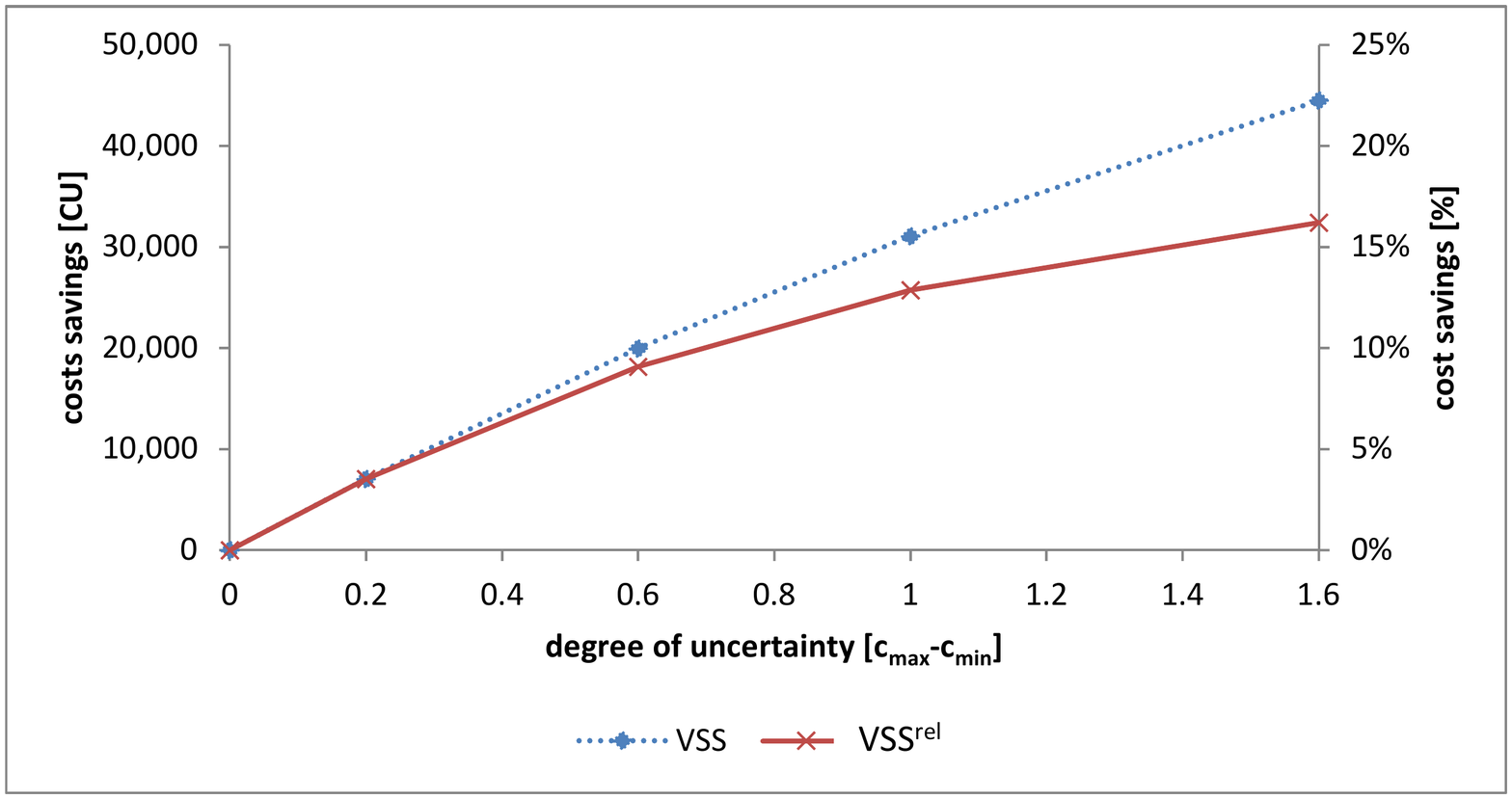}
\caption{Absolute and relative value of the stochastic solution for the instances of the basic instance structure, in dependence of the degree of uncertainty.}
\label{fig:VSS}       
\end{figure}

Whereas the last subsection investigated the difference between the predicted and the true costs of the EV solution, we turn now to the question of how much worse the EV solution is in comparison with the solution of the stochastic optimization problem. This latter difference indicates the value of taking the stochasticity of the demand into account instead of using the simplified deterministic EV model for the planning process. As the solution procedure for the stochastic optimization problem, we choose the MH approach.
As before, let $(z_{EV}^*,x_{EV}^*)$ and $E[G(z_{EV}^*,x_{EV}^*,\xi)]$ denote the EV solution and its expected external cost, respectively, and let $(z_{MH}^*,x_{MH}^*)$ and $E[G(z_{MH}^*,x_{MH}^*,\xi)])$ denote the MH solution and its expected external cost, respectively.
The difference\\
\noindent\begin{tabularx}{\linewidth}{@{} L {80mm}L {20mm}>{\raggedleft\arraybackslash}X@{}}\\
\multicolumn{2}{c}{
$\displaystyle VSS=E[G(z_{EV}^*,x_{EV}^*,\xi)] - E[G(z_{MH}^*,x_{MH}^*,\xi)]$
}&\tagarray\label{eq:EEVvsMH}\\
&
\end{tabularx}
is called the {\em value of the stochastic solution} (VSS); it describes the cost savings achieved by applying the stochastic solution approach instead of the deterministic one. We define the relative value of the stochastic solution as the quotient\\
\noindent\begin{tabularx}{\linewidth}{@{} L {80mm}L {20mm}>{\raggedleft\arraybackslash}X@{}}\\
\multicolumn{2}{c}{
$\displaystyle VSS^{rel}=(E[G(z_{EV}^*,x_{EV}^*,\xi)] - E[G(z_{MH}^*,x_{MH}^*,\xi)])/E[G(z_{EV}^*,x_{EV}^*,\xi)]$.
}&\tagarray\label{eq:EEVvsMHinprocent}\\
&
\end{tabularx}
In Figure \ref{fig:VSS}, the absolute and the relative VSS are depicted in dependence of the degree of uncertainty for the instances of the basic instance structure. Unsurprisingly, both the absolute and the relative VSS increase as uncertainty increases. The VSS rapidly grows with increasing degree of uncertainty, and it reaches about 16\%
for $[c_{min}, c_{max}]=[0.2,1.8]$. Applying a linear regression, we statistically confirm $(\alpha = 0.001)$ the intuitive conjecture that both the absolute and the relative VSS are positively correlated with the degree of uncertainty (for VSS and relative VSS, we get correlation coefficients of $R=0.936$ and $R=0.753$, respectively). In some other instance structures, the potential of the stochastic solution approach is even higher. For example, for large values of the number of projects and of the number of skills per resource, this gain reaches values between 30\% and 100\%.
The considerable size of the VSS confirms the need of applying stochastic optimization techniques in project scheduling and staffing.

\begin{figure}\centering
  \includegraphics[trim=25mm 45mm 25mm 42mm,,clip,scale=0.35]{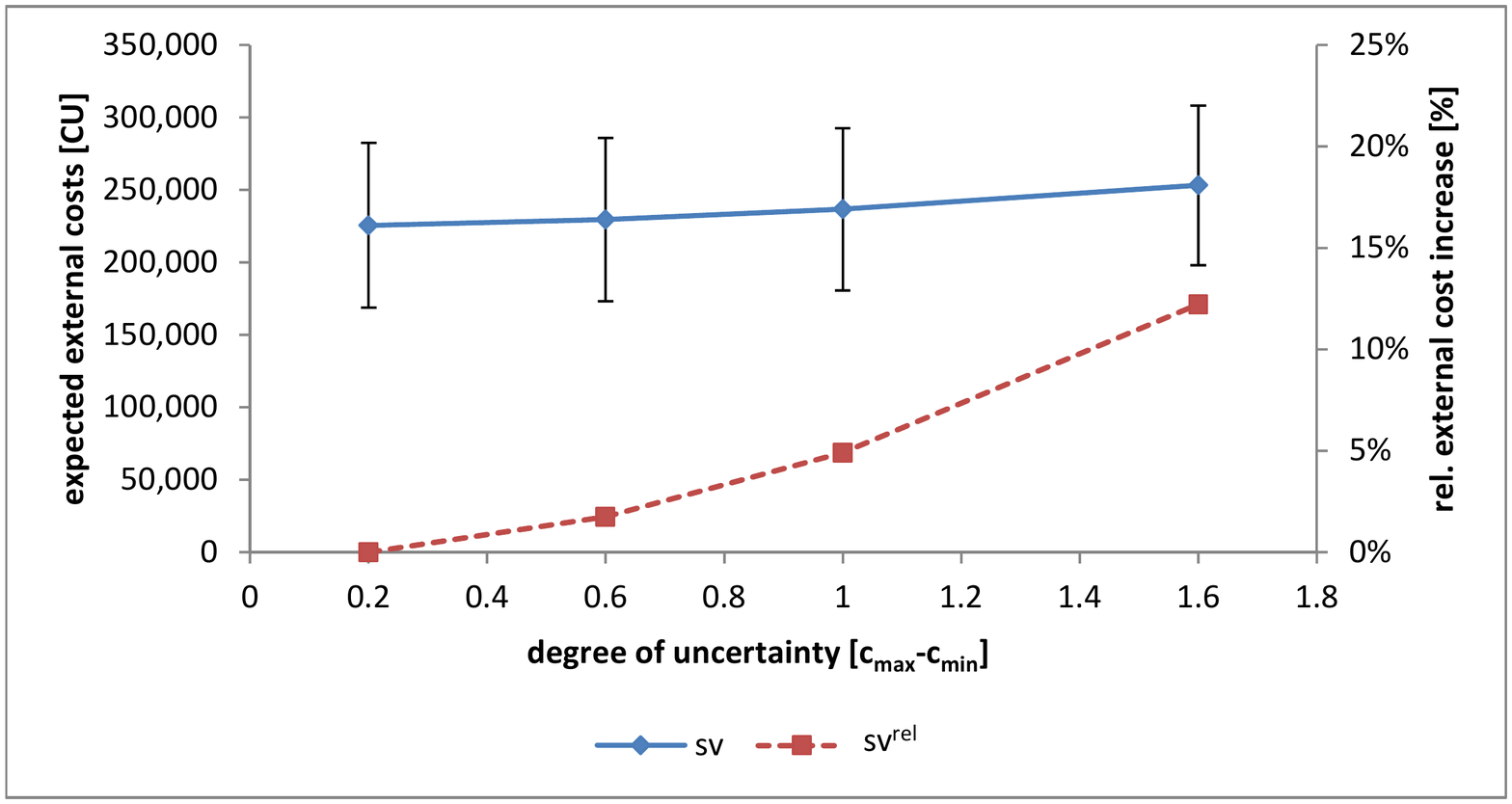}
\caption{Plot of (i) expected external costs and (ii) relative expected external cost increments compared to the situation with low stochasticity, $[c_{min},c_{max}]=[0.9,1.1]$, as a function of $c_{max} - c_{min}$, for the instances of the basic instance structure.}
\label{fig:degreeofuncertainty}       
\end{figure}

\subsection{Costs of Uncertainty}
\label{subsec:Resultsdegreeofuncertainty}

In Figure \ref{fig:degreeofuncertainty}, we plot (i) the expected external costs and (ii) the relative expected external cost increments, both in dependence of different degrees of uncertainty, for the instances of the basic scenario structure. In this Figure, the relative expected external cost increments indicate by which percentage the expected external costs increase if the situation $[c_{min},c_{max}]=[0.9,1.1]$ of low uncertainty is replaced by a higher uncertainty interval $[c_{min}, c_{max}]$ for the demand distributions.
As one would anticipate, the expected external costs increase as uncertainty increases. A linear regression for the absolute value of the expected external cost shows that the coefficient of correlation is $R=0.119$, with a p-value of 0.017 (significance at level $\alpha=0.05$).
Comparing the ``almost deterministic'' situation $[c_{min},c_{max}]=[0.9,1.1]$ to the situation $[c_{min},c_{max}]=[0.2,1.8]$ of poor information on demands, we find a gap of 12\%.
We would like to emphasize that a degree of uncertainty represented by $[c_{min},c_{max}]=[0.2,1.8]$ (i.e., a distribution allowing real efforts of work packages to exceed estimated efforts by up to 80\%) is not extreme from a applied point of view: In areas as software engineering or in the construction industry, even higher overruns occur.
Therefore, also the costs of uncertainty can be considerable in practice.

\subsection{Influence of parameters on external costs}
\label{subsec:ResultsInfluenceofthenumberofprojects}

In~\cite{Heimerl2010a}, the authors investigate in a deterministic context how parameters as the number of projects, the time window size etc.~influence the optimal costs. We shall extend now their results to the stochastic context of the present paper. In the present subsection, we use the {\em ceteris paribus design} explained in subsection~\ref{subsec:Testinstances}. That is, the parameters of the basic instance structure of Table~\ref{tab:basescenario} are applied, with the exception of modifying one single parameter among the parameters in Table~\ref{tab:experimentallevels}.

\begin{figure}
\hspace*{4ex}
{\includegraphics[trim=25mm 45mm 25mm 42mm,clip,scale=0.41]{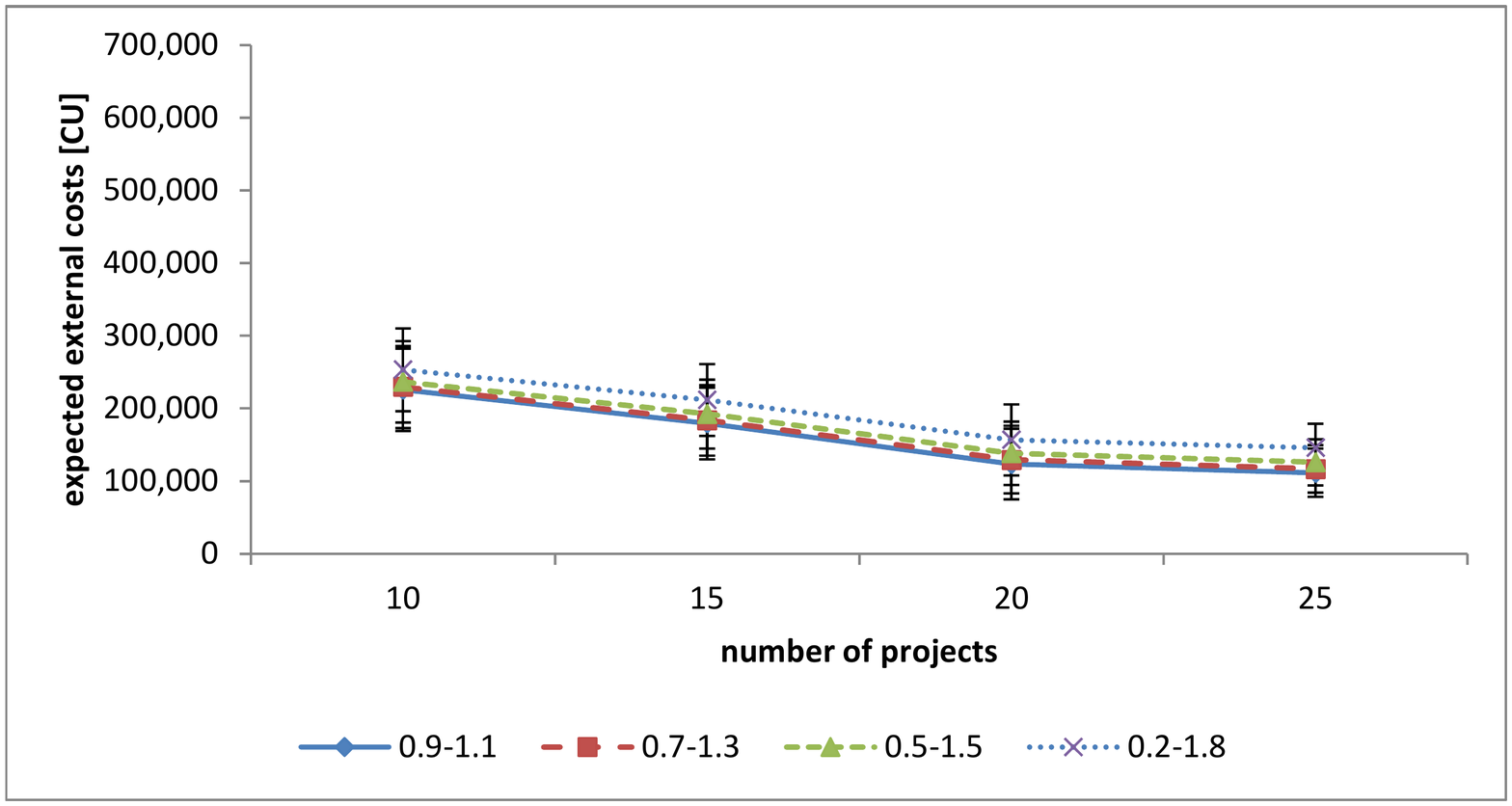}}
\caption{Plot of the expected external costs as a function of the number of projects $\vert P \vert$}
\label{fig:numberofproject}
\end{figure}
\begin{figure}
\hspace*{4ex}
{\includegraphics[trim=25mm 45mm 25mm 42mm,clip,scale=0.40]{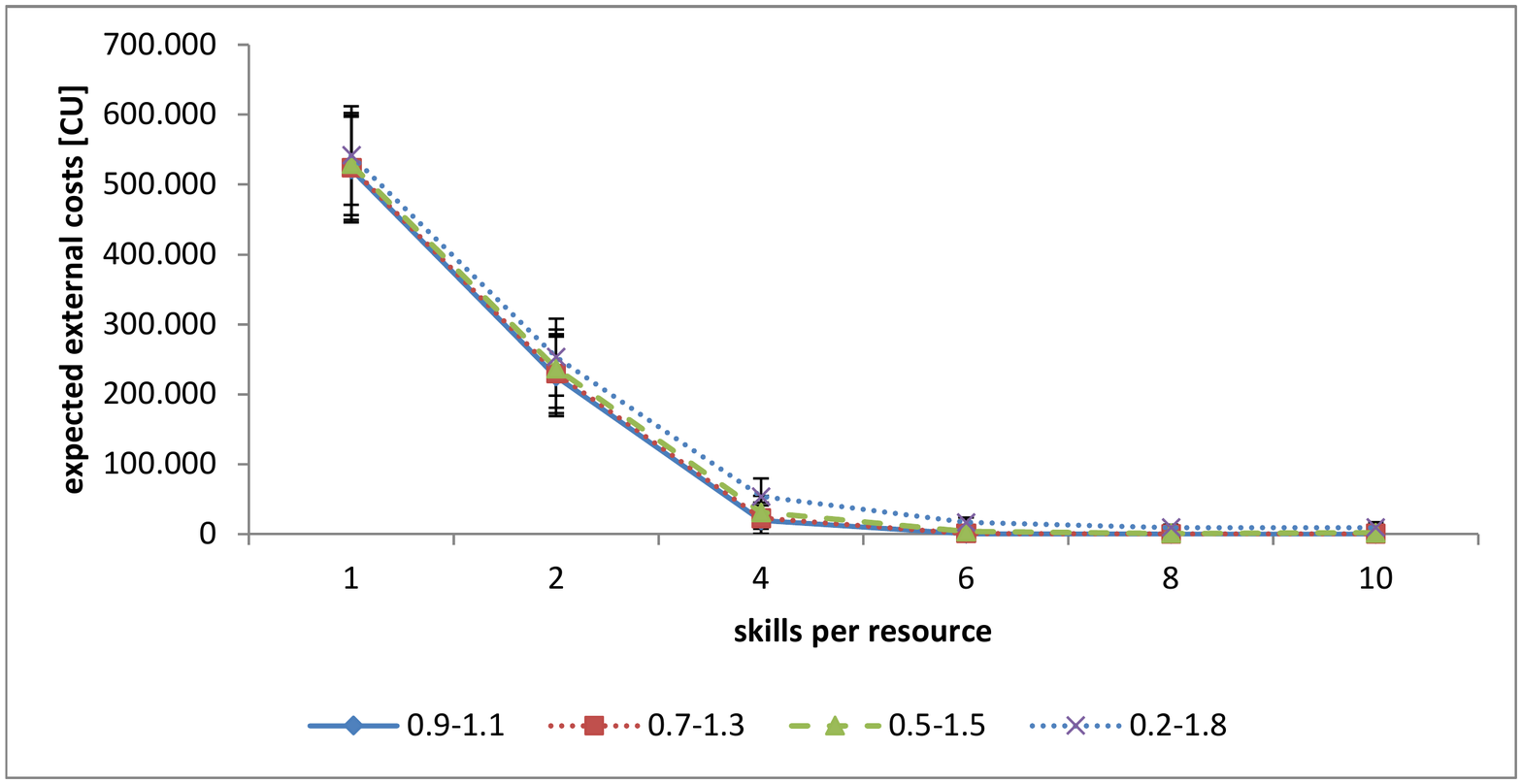}}
\caption{Plot of the expected external costs as a function of the number of skills per resource $\vert S_k \vert$}
\label{fig:skillsperresource}
\end{figure}

{\bf Influence of the number of projects.}
First, we investigate the influence of increasing the number of projects on the resulting expected external costs.
As~\cite{Heimerl2010a}, we keep the total resource demand constant while increasing the number of projects, which means that for a larger number of projects, the work packages become smaller. This change increases the flexibility of the planner, so one expects decreasing costs of the optimal solutions. This was confirmed indeed in \cite{Heimerl2010a} for the deterministic context.
We obtained similar results in the stochastic context: In Figure~\ref{fig:numberofproject}, we depict the expected external costs for four levels of uncertainty as a function of the number of projects (with fixed total demand). Additionally to the mean values, we plot the 95\% confidence interval to account for the randomness in test instance generation.
Apart from the already known effect that higher degree of uncertainty leads to higher external costs, one can see that all levels of uncertainty show the same behavior for an increasing number of projects:
The conjecture that a larger number of smaller work packages is easier to balance across the planning horizon than bigger and fewer work packages is confirmed.
Also as expected, we observe that from a certain value on,
the potential of the flexibility achieved by reducing work packages sizes diminishes.
{\bf Influence of the time window size.}
In a similar way as it was done for the influence of the number of projects, we investigated also the influence of the time window size on the costs. The results showed that in our test instances, an additional degree of freedom, i.e., a time window size change from zero to one, led to cost savings of around 20\%. A further cost reduction of up to 50\% was achieved by increasing the time window size~$\gamma$ from 1 to 2. For time window sizes larger than two time periods, the costs did not decrease further.
{\bf Influence of the number of skills per resource.}
In Figure \ref{fig:skillsperresource}, we plot the expected external costs in dependence of the number of skills per resource. Please note that the situation where the number of skills per resource is $\vert S_k \vert=10$ represents a situation when all resources posses all skills.
It can be observed that the external costs decrease monotonically with increasing $|S_k|$. Two special findings may be particularly important. First, in the context of the considered set of instances (basic instance structure, with modifications only with respect to $|S_k|$), a situation where all employees are extremely specialized, i.e., possess only one skill per person, leads to about the double expected external costs in comparison with a situation where the employees have two skills per person. A further investment in training that causes an increase in the number skills per person from two to four leads to another significant decrease in the expected external costs. Secondly, from a value of about four skills per person on, the investment in additional skills achieves only very limited cost savings. Of course, the quantitative amount of cost reduction depends on the specific characteristics of the test instance set (especially on the chosen parameter value for utilization, here: $\rho = 1$). Nevertheless, the results suggest the managerial insight that (i) multi-skilled resources can lead to significant cost savings, but (ii) companies should be aware that over-qualification allows no return on investment.
\vspace{1ex}

\begin{figure}
\hspace*{4ex}
{\includegraphics[trim=25mm 45mm 25mm 42mm,clip,scale=0.40]{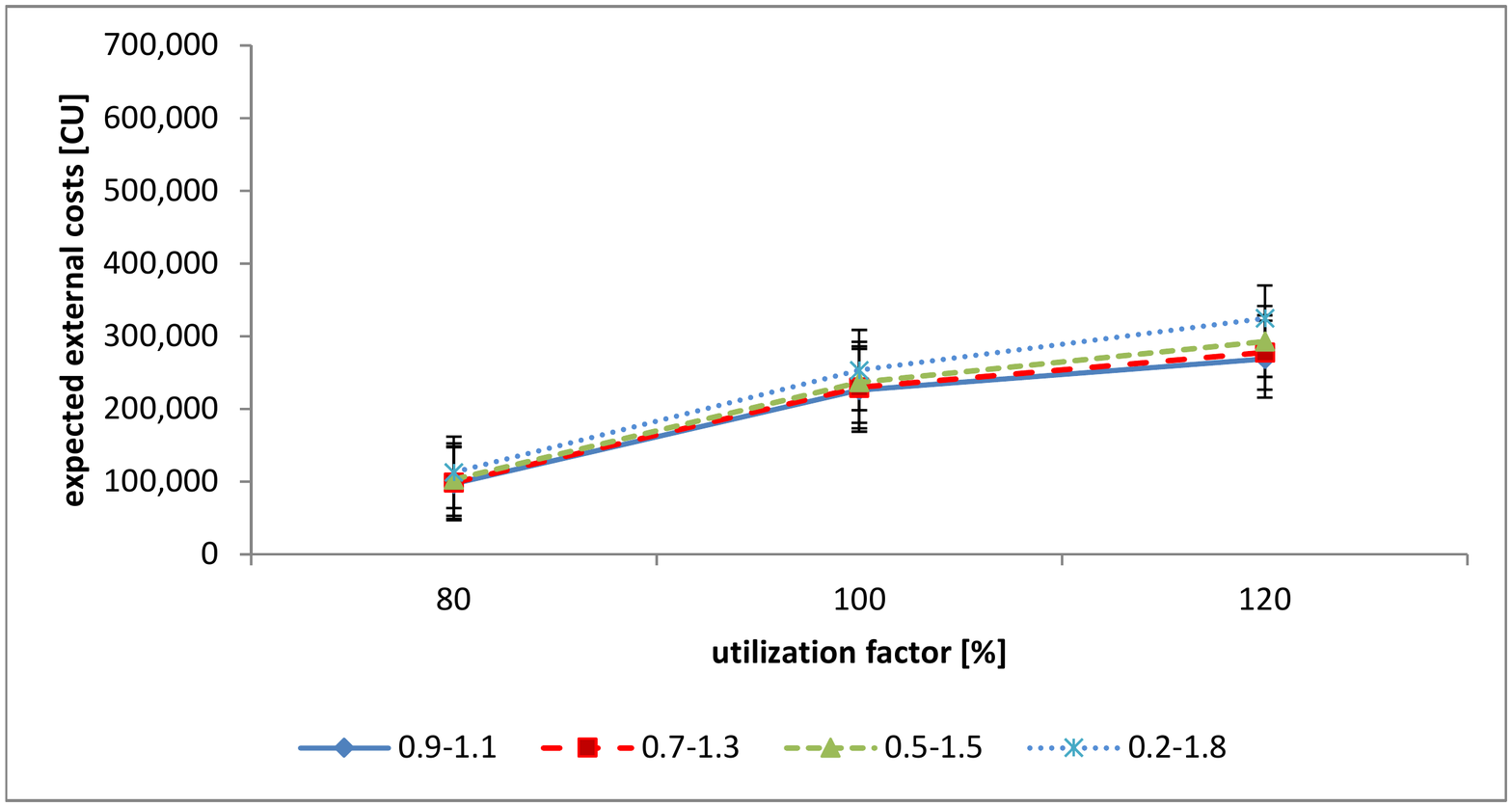}}
\caption{Plot of the expected external costs as a function of the utilization factor $\rho$}
\label{fig:workload}
\end{figure}

{\bf Influence of the utilization factor.}
Finally, in Figure~\ref{fig:workload}, the expected external costs are visualized as a function of the utilization factor~$\rho$, the
ratio of the expected demand to the available time capacities of the resources. As expected, Figure \ref{fig:workload} shows that the external costs increase as the utilization increases.
In more detail, we find that when starting at 100\% utilization, a 20\% decrease of utilization leads to a 57\% decrease of the expected external costs, whereas a 20\% increase of the utilization leads to a rise of costs by 23\%.


%
%
%
%

\subsection{Asymmetric work time distributions}

\begin{figure}\centering
 \includegraphics[trim=25mm 45mm 25mm 42mm,,clip,scale=0.35]{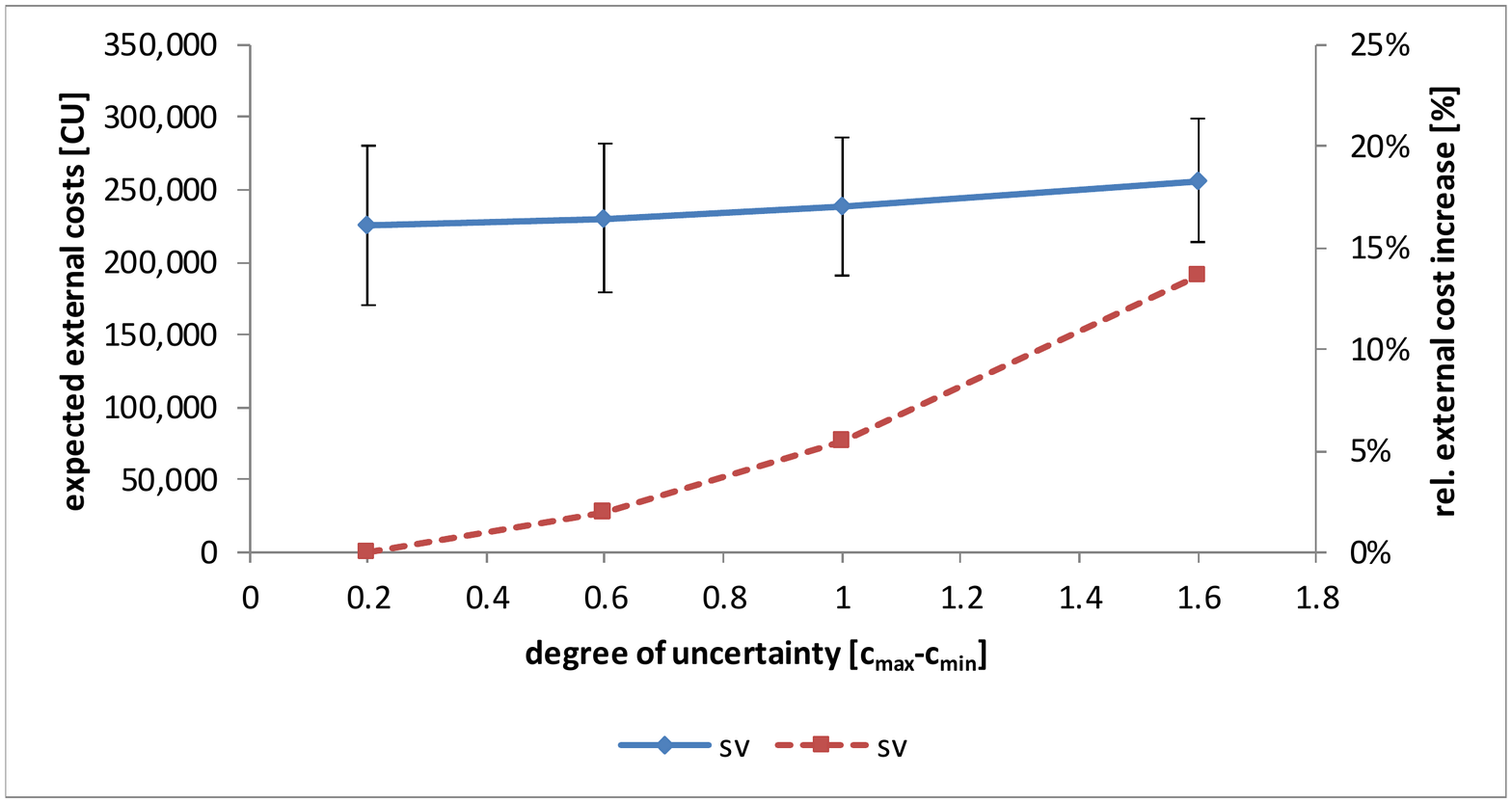}
\caption{Plot of (i) expected external costs and (ii) relative expected external cost increments compared to the situation with low stochasticity, $c_{max} - c_{min} = 0.2$, as a function of $c_{max} - c_{min}$, for the instances with right-skewed distribution.}
\label{fig:degreeofuncertaintyAS}       
\end{figure}

In the previous tests, we assumed that the triangular distribution of the variables $D_{psq}$ was {\em symmetric}. In practice, distributions of work times are often right-skewed. Therefore, we checked whether or not substantially different results were obtained by a replacement of the symmetric distributions by right-skewed triangular distributions. For this purpose, we took the basic instance structure and changed it as follows: The modal value was defined as $c_{mode}=c_{min}+(c_{max} - c_{min}) \cdot 0.25$, where $c_{min}$ and $c_{max}$ are the minimum and the maximum value of the distribution, respectively.
For a series of four test instances, we fixed the expected value $c_{expected}=( c_{min}+c_{max}+c_{mode})/3$ of the triangular distribution to 1. On this constraint, each of the four instances was constructed in such a way that the length $c_{max} - c_{min}$ of the support interval of the distribution was gradually increased, taking the values $0.2$, $0.6$, $1.0$ and $1.6$, respectively, which corresponds to increasing uncertainty.
Figure~\ref{fig:degreeofuncertaintyAS} shows a plot of the results. We see that the plot is very similar to Figure~\ref{fig:degreeofuncertainty}. We conclude that the skewedness of the distribution has little impact on the outcomes.

\section{Conclusion}
\label{sec:Conclusion and further research}

We developed two solution approaches for a stochastic project scheduling and staffing problem under uncertainty on required efforts. The scheduling decision is described by the choice of the start times of activities during the planning horizon. Each activity can consist of several work packages, where each work package requires a stochastic amount of effort exerted in one specific skill. The staff assignment decision matches available human resources (with heterogeneous skills) with the work package requirements of the activities. Demand for work time that is not covered by the assigned internal resources has to be satisfied by paying external work.

For our first solution approach, a ``matheuristic'' combination of a metaheuristic with a convex optimization procedure, we decompose the problem into a project scheduling subproblem and staffing subproblem. The project schedules are optimized by an iterated local search heuristic, using Variable Neighborhood Descent as a solution component. The search is guided to initially tackle time-periods leading to high external cost.
The staffing subproblem
is a convex optimization problem which we solve by the Frank-Wolfe algorithm. Different design variants of this algorithm are investigated in application to our problem.
Our second solution approach uses a Sample Average Approximation model of the stochastic optimization problem. This yields mixed-integer programming formulation that can be solved by CPLEX. For this second approach, a crucial parameter is the chosen sample size.

\begin{sloppypar}
Experimental results for synthetically generated test instances show that the matheuristic is a robust solution procedure that performs well for small as well as medium-sized test instances and provides solutions even in cases where the SAA model fails to return a feasible solution. Nevertheless, up to medium-sized instances, the sample average approximation is also a good choice. We find that deterministic planning replacing the stochastic optimization problem by the corresponding expected value problem bears the risk of drastically underestimating external costs.
Moreover, we demonstrate that the value of the stochastic solution, which is a measure for the cost savings achievable by a stochastic solution approach instead of using the deterministic expected value problem, is considerable, especially for moderate and high levels of uncertainty.
\end{sloppypar}

Our experiments confirmed some of the managerial insights that have been found (in a deterministic context) in \citet{Heimerl2010a} to be valid also in the context of uncertainty on required efforts. This holds especially for the potential as well as the limitations of multi-skilled resources, of small work packages and of large time windows with respect to possible cost savings. The main difference to the deterministic context of \citet{Heimerl2010a} is that a much more conservative planning strategy is necessary in presence of uncertainty in order to avoid excessive additional costs by required external work.

Finally, let us point out an important topic of future research: Whereas our approach only assumes substitutability of internal by external resources, the articles \citet{Ingels2017employee,ingels2017} focus on (within-skill or between-skill) substitutions of internal resources by other internal resources, which is another type of recourse action that is frequently deployed for absorbing work load peaks. It would be very interesting to extend the model presented here to these forms of reactive planning and to develop suitable solution techniques for such an extension.

Another topic of future research might be risk aversion. Our model is risk-neutral, which fits well to a situation of multi-project management with several small or medium-sized projects and a long-term perspective, but may be not appropriate anymore for a situation where very big problems cause specific risks. For such situations, extensions to risk-averse optimization approaches should be studied.
\vspace*{-2ex}

\section*{Acknowledgement}
This is a post-peer-review, pre-copyedit version of an article published in Journal of Scheduling. The final authenticated version is available online at \url{http://dx.doi.org/10.1007/s10951-018-0592-y}

\bibliographystyle{spbasic}      
\bibliography{bibfile}  

\end{document}